\documentclass[12pt,a4paper] {article}

\usepackage{amsthm,amssymb,amsmath,color,pdfsync,graphicx}
\usepackage[colorlinks,citecolor=blue, urlcolor=blue]{hyperref}
\usepackage{appendix}

\usepackage[utf8]{inputenc}
\usepackage[T1]{fontenc}
\usepackage{ae,aecompl}
\usepackage[paper=a4paper,textwidth=170mm]{geometry}

\newtheorem{rmk}{Remark}[section]

\newcommand{\be}{\begin{equation}}
\newcommand{\ee}{\end{equation}}
\newcommand{\bd}{\begin{displaymath}}
\newcommand{\ed}{\end{displaymath}}

\newcommand{\ba}{\begin{eqnarray}}
\newcommand{\ea}{\end{eqnarray}}
\newcommand{\ban}{\begin{eqnarray*}}
\newcommand{\ean}{\end{eqnarray*}}

\newcommand{\E} {I\!\! E}

\renewcommand{\Box}{\hfill\rule{0.25cm}{0.25cm}} 

\newtheorem{lemma}{Lemma}[section]
\newtheorem{theorem}{Theorem}[section]
\newtheorem{corollary}{Corollary}[section]

\def\1{{\bf 1}}

\newcounter{hypc}
\newcommand{\hyp}[1]{\stepcounter{hypc}\tag{$\mathbf{B_{\thehypc}}$}
\label{#1}}

\newcounter{rev}
\newcommand{\rev}[1]{\stepcounter{rev}\tag{$\mathbf{RE_{\therev}}$}
\label{#1}}

\newcounter{noisc}

\begin{document}

\title{Non-asymptotic Oracle Inequalities for the Lasso and Group Lasso in high dimensional logistic model}

%

\author{Marius Kw\'emou\footnote{Laboratoire Statistique et G\'enome UMR CNRS 8071- USC INRA, Universit\'e d'\'Evry Val d'Essonne, France; (2)LERSTAD, Universit\'e Gaston Berger de Saint-Louis, S\'en\'egal, e-mail: marius.kwemou@gmail.com}$^{~,1}$}

%
%

%
%

%
\maketitle
\begin{abstract}
We consider   the problem of estimating a function $f_{0}$ in logistic regression model. We propose to estimate this function $f_{0}$ by
a sparse approximation build as a linear combination of elements of a given dictionary of $p$ functions. This sparse approximation is selected by the Lasso or
Group Lasso procedure.
In this context, we state non asymptotic oracle inequalities for Lasso and Group Lasso under restricted eigenvalue assumption as introduced in \cite{BRT}.
\end{abstract}
\section*{Introduction}
During the last few years, logistic regression problems with more and more high-dimensional data occur in a wide
variety of scientific fields, especially in studies that attempt to find risk factors for disease and clinical outcomes. 
For example in gene expression data analysis or in genome wide association analysis the number $p$ of predictors may be
 of the same order or largely higher than the sample size $n$ (thousands $p$ of predictors  for  only a few dozens of individuals $n$, see
 for instance  \cite{manuel} or \cite{wu2009genome}).
In this context the considered model is often what we call here ``\textit{usual}" logistic regression. It is given by
\begin{equation}\label{usual}
\mathbb{P}(Y_{i}=1)=\pi(z_{i}^{T}\beta_{0})=\frac{\exp(z_{i}^{T}\beta_{0})}{1+\exp(z_{i}^{T}\beta_{0})},
\end{equation}
where one observes $n$ couples 
$(z_{1},Y_{1})$,\dots,$(z_{n},Y_{n})\in\mathbb{R}^{d}\times\{0,1\}$,  and $\beta_{0}$ is the unknown parameter to estimate.
Throughout the paper, we consider a fixed design setting (i.e $z_{1},\dots,z_{n}$ are considered deterministic).

In this paper, we consider  a more general logistic model described by
\begin{equation}\label{model}
 \mathbb{P}(Y_{i}=1)=\frac{\exp(f_{0}(z_{i}))}{1+\exp(f_{0}(z_{i}))},
\end{equation}
  where the outputs $Y_{i}\in\{0,1\},~i=1,\dots,n$ are independent and $f_{0}$ (not necessarily  linear) is  an unknown function (see \cite{non_par_logist}). 
 We aim at  estimating $f_{0}$ by constructing a suitable  approximation.
More precisely we estimate $f_{0}$ by a sparse approximation of linear combination of elements of a given dictionary of functions
$\mathbb{D}=\{\phi_{1},\dots,\phi_{p}\}$: $\hat{f}(.):=\sum_{j=1}^{p}\hat{\beta_{j}}\phi_{j}(.).$ 
Our purpose expresses the belief that, in many instances, even if $p$ is large, only a subset of $\mathbb{D}$ 
 may be needed to approximate $f_{0}$ well.
This construction can be done by minimizing the empirical risk.
However, it is well-known that with a large number of parameters in high dimensional data situations, direct minimization of empirical risk can lead to  
\textit{Overfitting}: the classifier can only behave well in training set, and can be bad in test set. The procedure would also be unstable: since empirical risk
is data dependent, hence random, small change in the data can lead to very different estimators.
 Penalization is used to overcome those drawbacks. One could use $\ell_{0}$ penalization, \textit{i.e.} penalized by the number of non zero coefficients 
(see for instance AIC, BIC \cite{akaike, schwarz}). Such a penalization would produce interpretable
models, but leads to non convex optimization and there is not efficient algorithm to solve this problem in high dimensional framework.
Tibshirani \cite{TR} proposes to use $\ell_{1}$ penalization, which is a regularization technique for simultaneous estimation and selection. This  penalization
 leads to convex optimization and is important from computational point of view (as well as from theoretical point 
of view).           
As a consequence of the optimality conditions, regularization by the $\ell_{1}$ penalty tends to produce some coefficients that
are exactly zero and shrink others, thus the name of Lasso (Least Absolute Shrinkage and Selection Operator). 
 There exist some algorithms to solve this convex problem,
 \textit{glmnet} (see \cite{Fried}), \textit{predictor-corector} (see \cite{park}) among the others.\\
 A related Lasso-type procedure is the Group Lasso, where the covariates are assumed to be
clustered in groups, and instead of $\ell_{1}$-penalty (summing the absolute values of each individual loading) the sum
of Euclidean norms of the loadings in each group is used. It shares the same kind of properties as the Lasso, but encourages
predictors to be selected in groups. This is useful when the set of predictors is partitioned into prescribed groups, only few being relevant
in the estimation process. Group Lasso has numerous applications : when categorical
predictors (factors) are present, the Lasso solution is not adequate since it only selects individual
dummy variables instead of whole factors. In this case, categorical variables are usually represented as groups of dummy variables. 
In speech and signal processing for example, the groups may represent different frequency bands~(see \cite{mcauley2005}).\\ 
\textbf{Previously known results.}
Recently, a great deal of attention has been focused on $\ell_{1}$-penalized based estimators. Most of this attention concerns regression models and $\ell_{1}$-penalized
least squares estimator of parameters in high dimensional linear and non linear additive regression. Among them one can cite 
\cite{bunea2006aggregation, bunea2007aggregation, bunea2007sparsity, masmen,bartlett2012, james2009dasso},
who have studied the Lasso for  linear model in nonparametric setting and proved sparsity oracle inequalities.
Similar sparsity oracle inequalities are proved in \cite{BRT},
 and those results hold under the so-called \textit{restricted eigenvalue assumption}
 on the Gram matrix. Those kind of results have been recently stated for the variants of the Lasso.
For instance Lounici et al. \cite{lounici2011} under a group version of \textit{restricted eigenvalue assumption}  stated oracle inequalities in 
linear gaussian noise model under Group sparsity. Those results lead to the refinements of their previous results for multi-task learning
(see \cite{lounici2009}). 
The behavior of the Lasso and Group Lasso regarding their  selection and estimation properties have been studied in :
 \cite{KF, meinshausen2006, zhao2006, osborne2000, zhang2008sparsity, meinshausen2009lasso} for Lasso in linear regression;
\cite{chesneau, nardi} for Group Lasso in linear regression; \cite{ravikumar, meier2009, huang2010} for additive models. 
Few results on the Lasso and Group Lasso concern logistic regression model. Most of them are asymptotic results and concern the ``usual" logistic 
regression model defined by (\ref{usual}). Zou \cite{zou2006} shows consistency in variable selection
for adaptive Lasso in generalized linear models when the number of covariables p is fixed. Huang et al. \cite{huang2008} prove sign consistency and estimation consistency for high-dimensional 
 logistic regression. Meir et al. \cite{mal} shown consistency for the Group Lasso in ``\textit{usual}"
logistic model (\ref{usual}).  To our knowledge there are only two non asymptotic results for the Lasso in logistic model : the first one is from
Bach \cite{bach10}, who  provided bounds for excess risk (generalization performance) and  estimation error in the case of ``\textit{usual}" logistic regression model
  under \textit{restricted eigenvalue assumption} on the weighted 
Gram matrix. The second one is from van de Geer \cite{van2008}, who established non asymptotic oracle inequality
 for Lasso in high dimensional generalized linear models with Lipschitz loss functions. Non asymptotic results concerning Group Lasso for logistic regression model
 have been established by  Negahban et al.  \cite{negahban}, with the assumption that $f_0$ is linear. 
 
 In this paper, we state  general non asymptotic oracle inequalities for the Lasso and Group Lasso in logistic model  within the
 framework of high-dimensional statistics.  We do not assume that $f_0$ is linear. We first state ``slow" oracle inequalities (see Theorem~\ref{theo02} and Theorem~\ref{theo0}) with no assumption on the Gram matrix, on the regressors nor on the margin.
Secondly we provide ``fast"  oracle inequalities (see Theorem~\ref{theo21} and Theorem~\ref{theo1}) under \textit{restricted eigenvalue assumption} 
and some technical assumptions on the regressors. In each case,
we give, as a consequence, the bounds for excess risk, $L_{2}(\frac{1}{n}\sum_{i=1}^{n}\delta_{z_{i}})$ and estimation errors for Lasso and Group Lasso in the ``usual" logistic regression.
Our non asymptotic results lead to an adaptive data-driven weighting of the $\ell_{1}$-norm (for the Lasso) and group norm (for the Group Lasso).
 
 This paper is organized as follows. In Section~\ref{s2}, we describe our weighted Group Lasso estimation procedure and state non asymptotic
 oracle inequalities for the Group Lasso estimator. In Section~\ref{s1} we describe our weighted Lasso estimation procedure and state non asymptotic
 oracle inequalities for the Lasso estimator. In Section~\ref{GLl} and Section~\ref{Ll}
we give as a consequence  the bounds for excess risk, $L_{2}(\frac{1}{n}\sum_{i=1}^{n}\delta_{z_{i}})$ and estimation errors for Lasso and Group Lasso in the ``\textit{usual}" logistic regression (\ref{usual})
. The proofs are gathered in Section~\ref{prouv} and Appendix.\\
\textbf{Definitions and notations}\\
Consider the matrix $X=\left(\phi_{j}(z_{i})\right)_{1\leq i \leq n,~ 1\leq j\leq p}$ and $\{G_{l},l=1,\dots,g\}$ the partition of $\{1,\dots,p\}$. 
 For any $\beta=(\beta_{1},\dots,\beta_{p})^{T}=(\beta^{1},\dots,\beta^{g})^{T}\in \mathbb{R}^{p},$ where $\beta^{l}=(\beta_{j})_{j\in G_{l}}$ for $l=1,\dots,g.$  
Let $f_{\beta}(.)=\sum_{j=1}^{p}\beta_{j}\phi_{j}(.)=\sum_{l=1}^{g}\sum_{j\in G_{l}}\beta_{j}\phi_{j}(.).$  
With our notations $$(f_{\beta}(z_{1}),\dots,f_{\beta}(z_{n}))^{T}=X\beta.$$ We define the group norm of $\beta$ as
$$ \lVert \beta \rVert_{2,q}=\left(\sum_{l=1}^{g}\left(\sum_{j\in G_{l}}\beta_{j}^{2}\right)^{\frac{q}{2}}\right)^{\frac{1}{q}}=
\left(\sum_{l=1}^{g}\lVert \beta^{l} \rVert_{2}^{q}\right)^{\frac{1}{q}}, $$
for every $1\leq q<\infty.$
 For $\beta\in \mathbb{R}^{p}$ $K(\beta)=\{j\in\{1,\dots,p\} : \beta_{j}\neq 0\}$ and $J(\beta)=\{l\in \{1,\dots,g\} :  \beta^{l}\neq 0 \},$ respectively the set of relevant coefficients
 (which characterizes the sparsity of the vector $\beta$)
 and the set of relevant groups. For all $\delta\in \mathbb{R}^{p}$ and a subset $I\subset \{1,\dots,p\},$ we denote by $\delta_{I}$
the vector in $\mathbb{R}^{p}$ that has the same coordinates as $\delta$ on $I$ and zero coordinates on the complement $I^{c}$ of $I.$ 
Moreover $|I|$ denotes the cardinality of $I$.
For all $h,f,g : \mathbb{R}^{d}\rightarrow R,$ we define the scalar products
$$\langle f,h\rangle_{n}=\frac{1}{n}\sum_{i=1}^{n}h(z_{i})f(z_{i}),$$
 and 
$$\langle f,h\rangle_{g}=\frac{1}{n}\sum_{i=1}^{n}h(z_{i})f(z_{i})\pi(g(z_{i}))(1-\pi(g(z_{i}))),~~\mbox{where}~~\pi(t)=\frac{\exp(t)}{1+\exp(t)}.$$
 We use the notation $$q_{f}(h)=\frac{1}{n}\sum_{i=1}^{n}h(z_{i})(Y_{i}-\pi(f(z_{i}))),$$ 
 $\lVert h\rVert_{\infty}=\max_{i}|h(z_{i})|$ and $\lVert h\rVert_{n}=\sqrt{\langle h,h\rangle_{n}}=\sqrt{\frac{1}{n}\sum_{i=1}^{n}h^{2}(z_{i})}$ 
which denote the $L_{2}(\frac{1}{n}\sum_{i=1}^{n}\delta_{z_{i}})$ norm (empirical norm). We consider empirical risk (logistic loss) for logistic model
\begin{equation}
\hat{R}(f)=\frac{1}{n}\sum_{i=1}^{n}\log(1+\exp(f(z_{i})))-Y_{i}f(z_{i}). 
\end{equation}
We denote by $R$ the expectation of $\hat{R}$ with respect to the distribution of $Y_{1},\dots,Y_{n},$ \textit{i.e}
$$ R(f)=\E(\hat{R}(f))=\frac{1}{n}\sum_{i=1}^{n}\log(1+\exp(f(z_{i})))-\E(Y_{i})f(z_{i}).$$
It is clear that $R(.)$ is a convex function and $f_{0}$ is a minimum of $R(.)$ when the model is well-specified (\textit{i.e.} when (\ref{model}) is satisfied).
Note that with our notations
\begin{equation}\label{equal}
 R(f)=\E(\hat{R}(f))=\hat{R}(f)+q_{f_{0}}(f).
\end{equation}
We shall use both the excess risk of  $f_{\hat{\beta}}$, $R(f_{\hat{\beta}})
-R(f_{0})$
and the prediction loss  $\lVert f_{\hat{\beta}}-f_{0}\rVert_{n}^{2}$ to evaluate the quality of the estimator.
 Note that  $R(f_{\hat{\beta}})$ corresponds
to the average Kullback-Leibler divergence to the best model when the
model is well-specified, and is common for the study of logistic regression. 

\section{Group Lasso for logistic regression model}\label{s2}
\subsection{Estimation procedure}
The goal is not to estimate the parameters of the  ``true" model (since there is no true parameter) but rather to construct an estimator that mimics the performance of the best model in a
given class, whether this model is true or not. Our aim is then  to estimate $f_{0}$ in Model (\ref{model}) by a linear combination of the functions of a dictionary 
$$\mathbb{D}=\{\phi_{1},\dots,\phi_{p}\},$$
where $\phi_{j} : \mathbb{R}^{d}\rightarrow \mathbb{R}$ and $p$ possibly  $>>n.$ The functions $\phi_{j}$ can be viewed as estimators of $f_{0}$ constructed from independent training
 sample, or estimators computed using $p$ different values of the tuning parameter of the same method. They can also be a collection of basis functions, that can 
approximate $f_{0},$ like wavelets, splines, kernels, etc... We implicitly assume that $f_{0}$ can be well approximated by a linear combination
$$f_{\beta}(.)=\sum_{j=1}^{p}\beta_{j}\phi_{j}(.),$$ 
 where $\beta$ has to be estimated.

In this section we assume that the set of relevant predictors
have  known group structure, for example in gene expression data these groups may be gene pathways,
or factor level indicators in categorical data. And we wish to achieves sparsity at the level of groups. This group sparsity assumption suggests us to use the Group Lasso
 method. We consider the Group Lasso for logistic regression (see~\cite{mal, yuan2006}), where predictors are included or excluded in groups.
The logistic Group Lasso is the minimizer of the following optimization problem
\begin{equation}\label{grl}
 f_{\hat{\beta}_{GL}} := \underset{f_{\beta}\in \varGamma}{\operatorname{argmin}}~\left\{\hat{R}(f_{\beta})+r\sum_{l=1}^{g}\omega_{l}\lVert \beta^{l} \rVert_{2}\right\},
\end{equation}
where
$$\varGamma \subseteq \left\{f_{\beta}(.)=\sum_{l=1}^{g}\sum_{j\in G_{l}}\beta_{j}\phi_{j}(.),~\beta\in \mathbb{R}^{p}\right\}.$$

  The tuning parameter $r>0$  is used to adjust the trade-off between minimizing the loss and finding a solution
which is sparse at the group level, i.e., to a vector $\beta$ such that $\beta^{l} = 0$ for
some of the groups $l\in \{1,\dots,g\}$. Sparsity is the consequence of the effect of non-differentiable penalty.
This penalty can be viewed as an intermediate between $\ell_{1}$ and $\ell_{2}$ type penalty, which has the attractive property 
that it does variables selection at the group level. The weights $\omega_{l}>0$, which we will define later, are used to control the amount of penalization per group.

\subsection{Oracle inequalities}
In this section we state non asymptotic oracle inequalities for excess risk and $L_{2}(\frac{1}{n}\sum_{i=1}^{n}\delta_{z_{i}})$ loss of Group Lasso estimator.
Consider the following assumptions :
\begin{align}
& \hyp{A1}  \mbox{There exists a constant}~  0 < c_{1} < \infty~ \mbox{such that } ~\max_{1\leq i \leq n}|f_{0}(z_{i})|\leq c_{1}.~~~~~~~~~~~~~~~~~~~~~~~~~~~~~~~~~~~\\
& \hyp{A4} \mbox{There exists a constant}~ 0 < c_{2} < \infty~\notag\mbox{such that} ~\max_{1\leq i \leq n}\max_{1\leq j \leq p}|\phi_{j}(z_{i})|\leq c_{2}.\\
&\hyp{A3} \mbox{There exists a  constant}~C_0 ~\mbox{such that the set}~~\varGamma=\varGamma(C_0)=\{f_\beta,\max_{1\leq i \leq n} |f_{\beta}(z_{i})|\leq C_{0}\} ~~\mbox{is non-empty}.~~~~~~~~~~~~~~~~~
\end{align}
Assumptions~(\ref{A1}) and (\ref{A3}) are technical assumptions useful to connect the excess risk and the $L_{2}(\frac{1}{n}\sum_{i=1}^{n}\delta_{z_{i}})$ loss 
(see Lemma~\ref{l81}). An assumption similar to (\ref{A1}) has been used in \cite{bunea2007aggregation} to prove oracle inequality in gaussian regression
model. The same kind of assumption 
 as (\ref{A3}) has been made in \cite{TS} to prove oracle inequality for support vector machine type with $\ell_{1}$ complexity
regularization. 
\begin{theorem}\label{theo02}
 Let  $f_{\hat{\beta}_{GL}}$ be the Group Lasso solution defined in (\ref{grl}) with $r\geq 1$ and 
\begin{equation}\label{omega2}
\omega_{l}=\frac{2\sqrt{|G_{l}|}}{n}\sqrt{\frac{1}{2} \underset{j\in G_{l}} \max \sum_{i=1}^{n}\phi_{j}^{2}(z_{i}) \left(x+\log{p}\right)}+\frac{2c_{2}\sqrt{|G_{l}|}}{3n}\left(x+\log{p}\right),
\end{equation}
where $x>0$. Under  Assumption~(\ref{A4}), with probability at least $1-2\exp(-x)$ we have
\begin{equation}\label{theo2eq1}
 R(f_{\hat{\beta}_{GL}})-R(f_{0})\leq \inf_{\beta\in \mathbb{R}^{p}}\left\{R(f_{\beta})-R(f_{0})+
2r\lVert \beta \rVert_{2,1}\underset{1\leq l \leq g}\max\omega_{l} \right\}.
\end{equation}
\end{theorem}
The first part of the right hand of Inequality (\ref{theo2eq1}) corresponds to the approximation error (bias). The selection of the dictionary can be very important to
minimize this approximation error. It is recommended to choose a dictionary $\mathbb{D}$ such that  $f_{0}$ could well be approximated by a linear combination of the 
functions of $\mathbb{D}.$  The
 second part of the right hand of Inequality (\ref{theo2eq1})  is the variance term and is usually referred as the rate of the oracle inequality.
In Theorem~\ref{theo02}, we speak about ``slow" oracle inequality, with the rate  at the order $\lVert \beta \rVert_{2,1}\sqrt{\log{p}/n}~~ \mbox{for any} ~~\beta.$ 
Moreover this is a sharp oracle inequality in the sense that there is a constant 1 in front of term 
$\underset{\beta\in \mathbb{R}^{p}}\inf\{R(f_{\beta})-R(f_{0})\}.$ This result is obtained without any assumption on the Gram matrix ($\Phi_{n}=X^{T}X/n$). 
In order to obtain oracle inequality with a  ``fast rate" of order  $\log{p}/n$ 
we need additional assumption on the restricted eigenvalue of the Gram matrix, namely the \textit{restricted eigenvalue assumption}.
\begin{align}
 &\rev{AA2} \mbox{For some integer}~ s ~\mbox{such that}~ 1\leq s \leq g~ \mbox{and a positive number} ~a_{0},~ \mbox{the following condition holds}\\
&\notag \mu_{1}(s,a_{0}) := \min_{K\subseteq \{1,...p\} : |K|\leq s}~\min_{\Delta\neq 0 : \lVert \Delta_{K^{c}} \rVert_{2,1}\leq a_{0}  \lVert \Delta_{K} \rVert_{2,1}}  
\frac{ \lVert X\Delta \rVert_{2}}{\sqrt{n} \lVert \Delta_{K} \rVert_{2}}>0 .
\end{align}
This is a natural extension to the Group Lasso of \textit{restricted eigenvalue assumption} introduced in~\cite{BRT} (or Assumption~(\ref{A2}) used below) for the usual Lasso. 
The only difference lies on the set where the minimum is taken : for the Lasso the minimum is taken over $\{\Delta\neq 0 : \lVert \Delta_{K^{c}} \rVert_{1}\leq a_{0}  \lVert \Delta_{K} \rVert_{1}\}$
 whereas for the Group Lasso the minimum is over
$\{\Delta\neq 0 : \lVert \Delta_{K^{c}} \rVert_{2,1}\leq a_{0}  \lVert \Delta_{K} \rVert_{2,1}\}.$ This assumption has already been used in \cite{lounici2009, lounici2011}
to prove oracle inequality for linear gaussian noise model under Group sparsity and for multi-task learning. To emphasize the dependency of
 Assumption~(\ref{AA2}) on $s$ and $a_{0}$ we will sometimes refer to it as $RE(s,a_{0}).$

\begin{theorem}\label{theo21}
 Let  $f_{\hat{\beta}_{GL}}$ be the Group Lasso solution defined in (\ref{grl}) with 
$\omega_{l}$ defined as in (\ref{omega2}).
 Fix $\eta>0$ and  $1\leq s \leq g$, assume  that (\ref{A1}), (\ref{A4}), (\ref{A3}) and (\ref{AA2}) are satisfied, with $a_{0}=3+4/\eta.$ 
 Thus with probability at least $1-2\exp(-x)$ we have
\begin{equation}\label{theo21eq1}
 R(f_{\hat{\beta}_{GL}})-R(f_{0})\leq (1+\eta)\inf_{f_{\beta}\in \varGamma}\left\{R(f_{\beta})-R(f_{0})+
\frac{c(\eta)|J(\beta)|r^{2}\left(\underset{1\leq l \leq g}\max\omega_{l}\right)^{2}}{c_{0}\epsilon_{0}\mu_{1}(s,a_{0})^{2}}   \right\},
\end{equation}
and
\begin{equation}\label{theo2eq2}
 \lVert f_{\hat{\beta}_{GL}}-f_{0}\rVert_{n}^{2}\leq \frac{c_{0}^{\prime}}{4c_{0}\epsilon_{0}}(1+\eta)\inf_{f_{\beta}\in \varGamma}\left\{\lVert f_{\beta}-f_{0}\rVert_{n}^{2}+
\frac{4c(\eta)|J(\beta)|r^{2}\left(\underset{1\leq l \leq g}\max\omega_{l}\right)^{2}}{c_{0}^{\prime}c_{0}\epsilon_{0}^{2}\mu_{1}(s,a_{})^{2}}   \right\}.
\end{equation}
Where $c(\eta)$ is a constant depending only on $\eta$; $c_{0}=c_{0}(C_{0},c_{1})$ and $c_{0}^{\prime}=c_{0}^{\prime}(C_{0},c_{1})$
 are constants depending on $C_{0}$ and $c_{1}$; $\epsilon_{0}=\epsilon_{0}(c_{1})$ is a constant depending on $c_{1}$; and $r\geq 1$ .
\end{theorem} 
 In Theorem~\ref{theo21}, the variance terms are of  order $\log{p}/n.$ Hence we say that the corresponding non asymptotic oracle inequalities have ``fast rates".
 For the best of our knowledge, Inequalities (\ref{theo2eq1}), (\ref{theo21eq1}) and (\ref{theo2eq2})  are the first non asymptotic oracle inequalities for the Group Lasso in logistic
regression model. These inequalities allow us to bound the prediction errors of Group Lasso by the best sparse approximation and a variance term. 
The major difference with existing results  concerning  Group Lasso for logistic regression model (see \cite{negahban,mal}) is that $f_0$ is not necessarily linear.

\begin{rmk}
  Our results remain true if we assume that we are in the ``neighborhood" of the target function. 
If we suppose that there exists $\zeta$ such that $\max_{1\leq i \leq n} |f_{\beta}(z_{i})-f_{0}(z_{i})|\leq \zeta$, then Lemma~\ref{l81} is still true.
\end{rmk}

\begin{rmk}\label{remark}
 The choice of the weights $\omega_{\ell}$ comes from Bernstein's inequality. We 
could also use the following weights 
$$\omega_{l}^{\prime}=\frac{2\sqrt{|G_{l}|}}{n}\sqrt{2 \underset{j\in G_{l}} \max \sum_{i=1}^{n}\mathbb{E}[\phi_{j}^{2}(z_{i})\epsilon_{i}^{2}] \left(x+\log{p}\right)}+
\frac{2\sqrt{|G_{l}|}\underset{1\leq i\leq n}\max\underset{j\in G_{l}}\max|\phi_j(z_i)|}{3n}\left(x+\log{p}\right),$$
where $\epsilon_i=Y_i-\mathbb{E}[Y_i]$, $i=1\dots n$.
Theorems~\ref{theo02} and \ref{theo21} still hold true with such weights $\omega_{l}^{\prime}.$
But these weights depend on the unknown function $f_{0}$ to be estimated through  $\E(\epsilon_{i}^{2})=\pi(f_{0}(z_{i})(1-\pi(f_{0}(z_{i})).$ 
This is the reason for using weights $\omega_{l}$ slightly greater than $\omega_{l}^{\prime}$. We also note that our weights are 
proportional to the square root of groups sizes, which is in acordance with the weights previously proposed for grouping strategies (see ~\cite{mal}).

\end{rmk}

\subsection{Special case : $f_{0}$ linear}\label{GLl}
In this section we assume that  $f_{0}$ is a linear function \textit{i.e.} $f_{0}(z_{i})=f_{\beta_{0}}(z_{i})=\sum_{l=1}^{g}\sum_{j\in G_{l}}\beta_{j}z_{ij}.$ 
Denote by $X=(z_{ij})_{1\leq i\leq n,1\leq j\leq p},$ the design matrix. Let $z_{i}=(z_{i1},\dots,z_{ip})^{T}$ be the ith row of the matrix $X$ and 
$z^{(j)}=(z_{1j},\dots,z_{nj})^{T}$ is jth column. For $i=1,\dots,n$
\begin{equation}\label{ligrp}
 \mathbb{P}(Y_{i}=1)=\frac{\exp(z_{i}^{T}\beta_{0})}{1+\exp(z_{i}^{T}\beta_{0})}.
\end{equation}
This corresponds to the ``\textit{usual}" logistic regression (\ref{usual}) \textit{i.e.} logistic model that allows linear dependency between $z_{i}$ and the distribution of $Y_{i}.$
In this context, the Group Lasso estimator of $\beta_{0}$ is defined by
\begin{equation}\label{gll}
\hat{\beta}_{GL} := \underset{\beta : ~f_{\beta}\in \varGamma}{\operatorname{argmin}}~\frac{1}{n}\sum_{i=1}^{n}\left\{\log(1+
\exp(z_{i}^{T}\beta))-Y_{i}z_{i}^{T}\beta\right\}+r\sum_{l=1}^{g}\omega_{l}\lVert \beta^{l} \rVert_{2}. 
\end{equation}
\begin{corollary}\label{col2}
 Let assumption \ref{AA2}(s,3) be satisfied and $|J(\beta_{0})|\leq s,$ where $1\leq s\leq g.$ Consider the Group Lasso estimator $f_{\hat{\beta}_{GL}}$ defined by (\ref{gll}) with
\begin{equation}\label{wl}
 \omega_{l}=\frac{2\sqrt{|G_{l}|}}{n}\sqrt{\frac{1}{2} \underset{j\in G_{l}} \max \sum_{i=1}^{n}z_{ij}^{2} \left(x+\log{p}\right)}+\frac{2c_{2}\sqrt{|G_{l}|}}{3n}\left(x+\log{p}\right)
\end{equation}
where $x>0.$ 
 Under the assumptions of Theorem~\ref{theo21},  with probability at least
$1-2\exp(-x)$ we have 
\begin{align}\label{col21}
 R(f_{\hat{\beta}_{GL}})-R(f_{\beta_{0}})&\leq \frac{9sr^{2}\left(\underset{1\leq l \leq g}\max\omega_{l}\right)^{2}}{\mu^{2}(s,3) c_{0}\epsilon_{0}}\\\label{col22}
\lVert f_{\hat{\beta}_{GL}}-f_{\beta_{0}}\rVert_{n}^{2}&\leq \frac{9sr^{2}\left(\underset{1\leq l \leq g}\max\omega_{l}\right)^{2}}{\mu^{2}(s,3) c_{0}^{2}\epsilon_{0}^{2}}\\\label{col23}
\lVert \hat{\beta}_{GL}-\beta_{0}\rVert_{2,1}&\leq  \frac{12rs\left(\underset{1\leq l \leq g}\max\omega_{l}\right)^{2} }{\mu^{2}(s,3) c_{0}\epsilon_{0}(\underset{1\leq l \leq g}\min\omega_{l})}\\\label{col24}
\lVert \hat{\beta}_{GL}-\beta_{0}\rVert_{2,q}^{q}&\leq  \left(\frac{12rs\left(\underset{1\leq l \leq g}\max\omega_{l}\right)^{2} }{\mu^{2}(s,3) c_{0}\epsilon_{0}(\underset{1\leq l \leq g}\min\omega_{l})}\right)^{q}~~~~\mbox{for all}~~ 1<q\leq 2.
\end{align}
\end{corollary}

\begin{rmk}\label{remlasso}
In logistic regression model~(\ref{logit}), if vector $\beta_{0}$  is sparse, \textit{i.e.} $|J(\beta_{0})|\leq s$, then Assumption~(\ref{AA2}) implies that $\beta_{0}$ is uniquely defined. Indeed, if there 
exists $\beta^{*}$ such that for $i=1,\dots,n$, $\pi(z_{i}^{T}\beta_{0})=\pi(z_{i}^{T}\beta^{*}),$ it follows that $X\beta_{0}=X\beta^{*}$ and $|J(\beta^{*})|\leq s.$ Then
according to assumption $RE(s,a_{0})$ with $a_{0}\geqslant 1,$ we necessarily have $\beta_{0}=\beta^{*}.$ Indeed  if $RE(s,a_{0})$ is satisfied with $a_{0}\geqslant 1,$
then $\min\{\lVert X\beta \rVert_{2} : |J(\beta)|\leq 2s, \beta\neq 0\}>0.$ 
\end{rmk} 

\begin{rmk} (\textbf{Theoretical advantage of Group Lasso over the Lasso})
 Concerning results on oracle inequality for the Group Lasso few results exist.
The first oracle inequality for the Group Lasso in the additive regression model
is due to~\cite{nardi}. Since then, some of these inequalities have been improved in Lounici et al. (2011) ~\cite{lounici2011}, concerning in particular the gain on order rate. More precisely,  Lounici et al. (2011)~\cite{lounici2011} have found a rate of order $\log{g}/n$ for Group Lasso in gaussian linear model, which is better than is corresponding  rate for the Lasso, $\log{p}/n$ (since $g\leq p$). This improvement seems mainly based on the assumption that the noise is gaussian. In our case (see proof of Theorem~\ref{theo02}, formula~(\ref{formula})) the empirical process involves non gaussian variables and thus their method should not apply in our context. 
However the probability that their results are true depends on $g$ whereas the probability that our results hold does not depend on $g$.

We can find the rate of order $\log{g}/n$  by  choosing this constant $x$ in the weights in a certain manner.
Indeed, let us assume (without loss of generality) that the groups are all of equal size $\lvert G_1\rvert=\dots=\lvert G_g\rvert = m$, so that $p =m.g$.
 Since the weights in (\ref{omega2}) are defined for  all $x>0$, if we take  $x= q\log{g}-\log{m}>0$  
where $q$ is a positive constant such that $g^{q}>m$. Then the weights in (\ref{omega2}) become
$$
\omega_{l}=\frac{2\sqrt{|G_{l}|}}{n}\sqrt{\frac{1}{2} \underset{j\in G_{l}} \max \sum_{i=1}^{n}\phi_{j}^{2}(z_{i}) \left[(1+q)\log{g}\right]}+\frac{2c_{2}\sqrt{|G_{l}}|}{3n}\left[(1+q)\log{g}\right],
$$
thus 
 $$\omega_{l}^{2} \sim  \frac{\log{g}}{n},$$
 and the results in Theorem ~\ref{theo02} and Theorem~\ref{theo21} hold with probability at least 
 $$1-2\frac{m}{g^{q}}.$$
 In the special case where the $g>2m$ these results are true for all $q>0$.
\end{rmk}
\subsection{Non bounded functions }
The results of Corollary~\ref{col2} are obtained (as the consequence of Theorem~\ref{theo21}) with the assumptions that $f_{\beta_0}$ 
and all $f_{\beta}\in\varGamma$ are bounded. In some situations these assumptions could not be verified.
In this section we will  establish the same results without assuming (\ref{A1}) or (\ref{A3}) \textit{i.e.} neither $f_{\beta_{0}}$ nor $f_{\beta}$ is bounded.
We consider the Group Lasso estimator defined in (\ref{gll}) and the following assumption :
\begin{align}
&\rev{AAA2} \mbox{For some integer}~ s ~\mbox{such that}~ 1\leq s \leq g~ \mbox{and a positive number}~a_{0},~ \mbox{the following condition holds}\\
&\notag\mu_{2}(s,a_{0}):= \min_{K\subseteq \{1,...p\} : |K|\leq s}~\min_{\Delta\neq 0 : \lVert \Delta_{K^{c}} \rVert_{2,1}\leq a_{0}  \lVert \Delta_{K} \rVert_{2,1}}  
\frac{ \Delta^{T}X^{T} DX\Delta}{ n \lVert \Delta_{K} \rVert_{2}^{2}}>0,\\
&\notag \mbox{where}~D=\mbox{Diag}\left(\mbox{var}(Y_{i})\right).
\end{align}
This is an extension of the Assumption~\ref{AA2} to the weighted Gram matrix $X^{T} DX/n.$
\begin{theorem}\label{theo22}
  Consider the Group Lasso estimator $f_{\hat{\beta}_{GL}}$ defined by (\ref{gll}) with $w_l$ defined as in (\ref{wl})
where $x>0.$  Set  $v=\underset{1\leq i\leq n}\max\underset{1\leq l \leq g}\max\lVert z_{i}^{l}\rVert_{2}.$
Let Assumptions~(\ref{A4}) and (\ref{AAA2}) be satisfied with $$a_{0}=\frac{3\underset{1\leq l \leq g}\max\omega_{l}}{\underset{1\leq l \leq g}\min\omega_{l}}.$$  If 
 $r(1+a_0)^{2}\underset{1\leq l \leq g}\max\omega_{l}\leq \frac{\mu_{2}^{2}}{3v|J|},$
with probability at least $1-2\exp(-x)$ we have 
\begin{align}\label{theo221}
 R(f_{\hat{\beta}_{GL}})-R(f_{\beta_{0}})&\leq \frac{9(1+a_0)^{2}J(\beta_{0})|r^{2}\left(\underset{1\leq l \leq g}\max\omega_{l}\right)^{2}}{\mu_{2}^{2}(s,3)}\\ \label{theo222}
\lVert \hat{\beta}_{GL}-\beta_{0}\rVert_{2,1}&\leq  \frac{6(1+a_0)^{2} |J(\beta_0)|r\left(\underset{1\leq l \leq g}\max\omega_{l}\right)}{\mu_{2}^{2}(s,3)}\\\label{theo223}
\lVert \hat{\beta}_{GL}-\beta_{0}\rVert_{2,q}^{q}&\leq  \left( \frac{6(1+a_0)^{2} |J(\beta_0)|r\left(\underset{1\leq l \leq g}\max\omega_{l}\right)}{\mu_{2}^{2}(s,3)}\right)^{q}~~~~\mbox{for all}~~  1<q\leq 2.
\end{align}
Moreover if we assume that there exists  $0<\epsilon_{0}\leq 1/2$ such that $$\epsilon_{0}\leq \pi(f_{\beta_{0}}(z_{i}))[1-\pi(f_{\beta_{0}}(z_{i}))]~~~~~\mbox{for all}~~  i=1,\dots,n $$ then,
\begin{equation}\label{theo224}
\lVert X\hat{\beta}_{GL}-X\beta_{0}\rVert_{n}^{2}\leq \frac{36(1+a_0)^{2}|J(\beta_{0})|r^{2}\left(\underset{1\leq l \leq g}\max\omega_{l}\right)^{2}}{\mu^{2}(s,3)\epsilon_{0}} .
\end{equation}
\end{theorem}
Inequalities~(\ref{theo222}) and (\ref{theo223}) are the extensions of the results in~\cite{bach10} for the Lasso to Group Lasso in logistic regression model.

In this section we studied some properties of the Group Lasso. However the  Group Lasso is based on prior knowledge that
the set of relevant predictors have  known group structure. If this group sparsity condition is not satisfied, the sparsity can be achieve by simply using the Lasso.  
We will show in the next section how to adapt the results of this section to the Lasso.
\section{Lasso for logistic regression}\label{s1}
\subsection{Estimation procedure}
The Lasso estimator $f_{\hat{\beta}_{L}}$ is defined as a minimizer of the following $\ell_{1}$-penalized empirical risk 
\begin{equation}\label{Lasso}
 f_{\hat{\beta}_{L}} := \underset{f_{\beta}\in \varGamma}{\operatorname{argmin}}~\left\{\hat{R}(f_{\beta})+r\sum_{j=1}^{p}\omega_{j}|\beta_{j}|\right\},
\end{equation}
where the minimum is taken over the set 
$$\varGamma  \subseteq \left\{f_{\beta}(.)=\sum_{j=1}^{p}\beta_{j}\phi_{j}(.),~\beta=(\beta_{1},\dots,\beta_{p})\in \mathbb{R}^{p}\right\}$$
and $\omega_{j}$ are positive weights to be specified later. The ``classical" Lasso penalization corresponds to $\omega_{j}=1,$ where $r$ is
the tuning parameter which makes balance between goodness-of-fit and sparsity. The Lasso estimator has the property that it does predictors 
selection and estimation at the same time. Indeed
 for large values of $\omega_{j}$, the related components $\hat{\beta}_{j}$ are set exactly to $0$ and the other are shrunken toward zero.

\subsection{Oracle inequalities}\label{slasso}
In this section we provide non asymptotic oracle inequalities for the Lasso in logistic regression model.
\begin{theorem}\label{theo0}
 Let $f_{\hat{\beta}_{L}}$ be the $\ell_{1}$-penalized minimum defined in (\ref{Lasso}). Let Assumption~(\ref{A4}) be satisfied. 
\begin{enumerate}
\item[A-)] Let $x>0$ be fixed and $r\geq 1$. For $j=\{1,\dots,p\},$ let
\begin{equation}\label{omega}
\omega_{j}=\frac{2}{n}\sqrt{\frac{1}{2}\sum_{i=1}^{n}\phi_{j}^{2}(z_{i})(x+\log{p})}+\frac{2c_{2}(x+\log{p})}{3n}.
\end{equation}
Thus with probability at least $1-2\exp(-x)$ we have
\begin{equation*}
 R(f_{\hat{\beta}_{L}})-R(f_{0})\leq \inf_{\beta\in \mathbb{R}^{p}}\left\{R(f_{\beta})-R(f_{0})+
   2\lVert \beta \rVert_{1}r\underset{1\leq j \leq p}\max\omega_{j}\right\}.
\end{equation*}
 \item[B-)] Let $A>2\sqrt{c_{2}}$. For $j=\{1,\dots,p\},$ let $\omega_{j}=1,$ and
 $$r=A\sqrt{\frac{\log{p}}{n}}.$$
  Thus with probability at least $1-2p^{1-A^{2}/4c_{2}}$ we have
\begin{equation*}
 R(f_{\hat{\beta}_{L}})-R(f_{0})\leq \inf_{\beta\in \mathbb{R}^{p}}\left\{R(f_{\beta})-R(f_{0})+
2A\lVert \beta \rVert_{1}r\sqrt{\frac{\log{p}}{n}}   \right\}.
\end{equation*}
\end{enumerate}
\end{theorem}

As previously, the variance terms are of order $\lVert \beta \rVert_{1}\sqrt{\log{p}/n}$ for any $\beta.$ Hence these are sharp oracle inequalities with ``slow" rates.  
These results are obtained without any assumption on the Gram matrix. To obtain oracle inequalities with a ``fast rate", of order $\log{p}/n$, we need the 
restricted eigenvalue condition. 
\begin{align}
&\rev{A2}  \mbox{For some integer} ~s~ \mbox{such that}~ 1\leq s \leq p ~\mbox{and a positive number}~  a_{0},
\mbox{the following condition holds}\\
&\notag \mu(s,a_{0}) := \min_{K\subseteq \{1,...p\} :|K|\leq s}~\min_{\Delta\neq 0 : \lVert \Delta_{K^{c}} \rVert_{1}\leq a_{0}  
\lVert \Delta_{K} \rVert_{1}} \frac{ \lVert X\Delta \rVert_{2}}{\sqrt{n} \lVert \Delta_{K} \rVert_{2}}>0 .
\end{align}

This assumption has been introduced in~\cite{BRT}, where several sufficient conditions
for this assumption are described. This condition is known to be one of the weakest to derive ``fast rates" for the Lasso. For instance conditions on 
the Gram matrix used to prove oracle inequality in \cite{bunea2006aggregation, bunea2007aggregation, bunea2007sparsity} are more 
restrictive than \textit{restricted eigenvalue assumption}. 
In those papers either $\Phi_{n}$ is positive definite, or mutual coherence condition 
is imposed. We refer to \cite{van2009} for a complete comparison of the assumptions used to prove oracle inequality for the Lasso.  
 Especially it is proved that \textit{restricted eigenvalue assumption} is weaker than the neighborhood stability or  irrepresentable  condition. 
\begin{theorem}\label{theo1}
 Let $f_{\hat{\beta}_{L}}$ be the $\ell_{1}$-penalized minimum defined in (\ref{Lasso}). 
Fix $\eta>0$ and $1\leq s \leq p$.
Assume that (\ref{A1}), (\ref{A4}), (\ref{A3}) and (\ref{A2}) are satisfied, with $a_{0}=3+4/\eta$.
\begin{enumerate}
\item[A-)] Let $x>0$ be fixed and $r\geq 1$. For $j=\{1,\dots,p\},$ $\omega_{j}$ defined as in (\ref{omega}).
Thus with probability at least $1-2\exp(-x)$ we have
\begin{equation}\label{theo1eq1}
 R(f_{\hat{\beta}_{L}})-R(f_{0})\leq (1+\eta)\inf_{f_{\beta}\in \varGamma}\left\{R(f_{\beta})-R(f_{0})+
\frac{c(\eta)|K(\beta)|r^{2}\left(\underset{1\leq j \leq p}\max\omega_{j}\right)^{2}}{c_{0}\epsilon_{0}\mu^{2}(s,3+4/\eta)}   \right\},
\end{equation}
and
\begin{equation}\label{theo1eq2}
 \lVert f_{\hat{\beta}_{L}}-f_{0}\rVert_{n}^{2}\leq \frac{c_{0}^{\prime}}{4c_{0}\epsilon_{0}}(1+\eta)\inf_{f_{\beta}\in \varGamma}\left\{\lVert f_{\beta}-f_{0}\rVert_{n}^{2}+
\frac{4c(\eta)|K(\beta)|r^{2}\left(\underset{1\leq j \leq p}\max\omega_{j}\right)^{2}}{c_{0}^{\prime}c_{0}\epsilon_{0}^{2}\mu^{2}(s,3+4/\eta)}   \right\}.
\end{equation}
 \item[B-)] Let $A>2\sqrt{c_{2}}$. For $j=\{1,\dots,p\},$ let $\omega_{j}=1,$ and
 $$r=A\sqrt{\frac{\log{p}}{n}}.$$
  Thus with probability at least $1-2p^{1-A^{2}/4c_{2}}$ we have
\begin{equation}\label{theo1eq3}
 R(f_{\hat{\beta}_{L}})-R(f_{0})\leq (1+\eta)\inf_{f_{\beta}\in \varGamma}\left\{R(f_{\beta})-R(f_{0})+
\frac{A^{2}c(\eta)}{c_{0}\epsilon_{0}\mu^{2}(s,3+4/\eta)}\frac{|K(\beta)|r^{2}\log{p}}{n}   \right\},
\end{equation}
and
\begin{equation}\label{theo1eq4}
 \lVert f_{\hat{\beta}_{L}}-f_{0}\rVert_{n}^{2}\leq \frac{c_{0}^{\prime}}{4c_{0}\epsilon_{0}}(1+\eta)\inf_{f_{\beta}\in \varGamma}\left\{\lVert f_{\beta}-f_{0}\rVert_{n}^{2}+
\frac{4c(\eta)A^{2}}{c_{0}^{\prime}c_{0}\epsilon_{0}^{2}\mu^{2}(s,3+4/\eta)}\frac{|K(\beta)|r^{2}\log{p}}{n}    \right\}.
\end{equation}
\end{enumerate}

In both cases $c(\eta)$ is a constant depending only on $\eta$; $c_{0}=c_{0}(C_{0},c_{1})$ and $c_{0}^{\prime}=c_{0}^{\prime}(C_{0},c_{1})$
 are constants depending on $C_{0}$ and $c_{1}$; and $\epsilon_{0}=\epsilon_{0}(c_{1})$ is a constant depending on $c_{1}$.
\end{theorem}
In this theorem the variance terms are of order $|K(\beta)|\log{p}/n.$  Such order in sparse oracle inequalities usually refer to ``fast rate". This rate
is of same kind of the one obtain in \cite{BRT} for linear regression model.
For the best of our knowledge, (\ref{theo1eq2}) and (\ref{theo1eq4}) are the first non asymptotic oracle inequalities for the $L_{2}(\frac{1}{n}\sum_{i}^{n}\delta_{z_i})$ norm in logistic model. 
Some non asymptotic oracle inequalities for excess risk like (\ref{theo1eq1}) or (\ref{theo1eq3}) have been established in~\cite{van2008} under different assumptions.
Indeed, she stated oracle inequality for high dimensional generalized linear model with Lipschitz loss function, where logistic regression is 
a particular case. Her result assumes to be hold in the ``neighborhood" of the target function, while our result is true for all
bounded functions. Note also that our results hold under $RE$ condition, which can be seen as empirical version of
Assumption C in~\cite{van2008}. 
The confidence (probability that result holds true) of Inequality~(\ref{theo1eq1}) does not depend on $n$ or $p$ while the confidence of her results depends on $n$ and $p$.
Moreover, the weights we proposed from Bernstein's inequality  are different
and easy to interpret. 

%

\subsection{Special case : $f_0$ linear}\label{Ll}
In this section we assume that  $f_{0}$ is a linear function that is $f_{0}(z_{i})=f_{\beta_{0}}(z_{i})=\sum_{j=1}^{p}\beta_{0j}z_{ij}=z_{i}^{T}\beta_{0},$ 
where $z_{i}=(z_{i1},\dots,z_{ip})^{T}.$ Denote $X=(z_{ij})_{1\leq i\leq n, 1\leq j\leq p}$ the design matrix. Thus
for $i=1,\dots,n$
\begin{equation}\label{logit}
 \mathbb{P}(Y_{i}=1)=\pi(z_{i}^{T}\beta_{0})=\frac{\exp(z_{i}^{T}\beta_{0})}{1+\exp(z_{i}^{T}\beta_{0})}.
\end{equation}
The Lasso estimator of $\beta_{0}$ is thus defined as
\begin{equation}\label{l}
\hat{\beta}_{L} := \underset{\beta : ~f_{\beta}\in \varGamma}{\operatorname{argmin}}~\left\{\frac{1}{n}\sum_{i=1}^{n}\left\{\log(1+\exp(z_{i}^{T}\beta))-Y_{i}z_{i}^{T}\beta\right\}+
r\sum_{j=1}^{p}\omega_{j}|\beta_{j}|\right\}. 
\end{equation}
 When the design matrix $X$ has full rank, the solution of optimization Problem~(\ref{l})
is usually unique.  When $p>>n$ this infimum might not be unique.

\begin{corollary}\label{col1}
 Let assumption RE(s,3) be satisfied and $|K(\beta_{0})|\leq s,$ where $1\leq s\leq p$. Consider the Lasso estimator $f_{\hat{\beta}_{L}}$
 defined by (\ref{l}) with  
$$\omega_{j}=\frac{2}{n}\sqrt{\frac{1}{2}\sum_{i=1}^{n}z_{ij}^{2}(x+\log{p})}+\frac{2c_{2}(x+\log{p})}{3n}$$
 Under the assumptions of Theorem~\ref{theo1} with probability at least $1-\exp(-x)$ we have 
\begin{align}\label{col1eq1}
 R(f_{\hat{\beta}_{L}})-R(f_{\beta_{0}})&\leq \frac{9sr^{2}\left(\underset{1\leq j \leq p}\max\omega_{j}\right)^{2}}{\mu^{2}(s,3) c_{0}\epsilon_{0}}\\\label{col1eq2}
\lVert f_{\hat{\beta}_{L}}-f_{\beta_{0}}\rVert_{n}^{2}&\leq \frac{9s^{2}r^{2}\left(\underset{1\leq j \leq p}\max\omega_{j}\right)^{2}}{\mu^{2}(s,3) c_{0}^{2}\epsilon_{0}^{2}}\\\label{col1eq3}
\lVert \hat{\beta}_{L}-\beta_{0}\rVert_{1}&\leq  \frac{12 sr \left(\underset{1\leq j \leq p}\max\omega_{j}\right)^{2}}{\mu^{2}(s,3) c_{0}\epsilon_{0}\left(\underset{1\leq j \leq p}\min\omega_{j}\right)}\\\label{col1eq4}
\lVert \hat{\beta}_{L}-\beta_{0}\rVert_{q}^{q}&\leq  \left(\frac{12 sr \left(\underset{1\leq j \leq p}\max\omega_{j}\right)^{2}}{\mu^{2}(s,3) c_{0}\epsilon_{0}\left(\underset{1\leq j \leq p}\min\omega_{j}\right)}\right)^{q}~~~~\mbox{for all}~~ 1<q\leq 2.
\end{align}
If $r=A\sqrt{\log{p}/n}$ and $\omega_{j}=1$ for all $j\in \{1,\dots,p\}$ we have the same results with probability at least
$1-2p^{1-A^{2}/4c_{2}}.$
\end{corollary}
Line~(\ref{col1eq1}) and Line~(\ref{col1eq3}) of the corollary  are similar to those of Theorem~5 in \cite{bach10}. Note that, up to differences in constant factors, the rates obtained in this corollary are the same as those obtained
 in Theorem 7.2 in \cite{BRT} for linear model with an s-sparse vector. Remark~\ref{remlasso} remains true in this section.

\section{Conclusion}
In this paper we stated non asymptotic oracle inequalities for the Lasso and Group Lasso.
Our results are non asymptotic : the number $n$ of observations is fixed while the number $p$ of covariates can grow with respect to $n$ and can be much larger
than $n$.  The major difference with existing results concerning Group Lasso or Lasso for logistic regression model is that we do not assume that $f_0$ is linear.  First we provided sharp oracle inequalities for excess risk, with ``slow" rates, with no assumption on the Gram matrix, on the regressors nor on the margin. Secondly, under RE condition we provided
 ``fast" oracle inequalities for excess risk and $L_{2}(\frac{1}{n}\sum_{i=1}^{n}\delta_{z_{i}})$ loss. We also provided as a consequence of oracle inequalities the bounds for excess
risk, $L_{2}(\frac{1}{n}\sum_{i=1}^{n}\delta_{z_{i}})$ error and estimation error in the case where the true function $f_{0}$ is linear (``\textit{usual}" logistic regression (\ref{usual})).
\section*{Acknowledgements}
We would like to thank Marie-Luce Taupin  for the careful reading of the manuscript and for her helpful comments. We also thank Sarah Lemler for helpful discussions.
\section{Proofs of main results}\label{prouv}
\subsection{Proof of Theorem~\ref{theo02}}
Since $\hat{\beta}_{GL}$ is the minimizer of $\hat{R}(f_{\beta})+r\sum_{l=1}^{g}\omega_{l}\lVert \beta^{l} \rVert_{2}$, we get
$$R(f_{\hat{\beta}_{GL}})-\frac{1}{n}\varepsilon^{T}X\hat{\beta}_{GL} +r\sum_{l=1}^{g}\omega_{l}\lVert \hat{\beta}^{l}_{GL} \rVert_{2}\leq R(f_{\beta})-\frac{1}{n}\varepsilon^{T}X\beta +
r\sum_{l=1}^{g}\omega_{l}\lVert \beta^{l} \rVert_{2},$$
where $\epsilon =(\epsilon_1,\dots,\epsilon_n)^{T}$ with $\epsilon_i = Y_i - \mathbb{E}[Y_i]$ for $i=1,\dots,n$.
By applying Cauchy-Schwarz inequality, we obtain 
\begin{eqnarray}
 \notag R(f_{\hat{\beta}_{GL}})-R(f_{0})&\leq & R(f_{\beta})-R(f_{0})+
\sum_{l=1}^{g}\frac{1}{n}\sqrt{\sum_{j\in G_{l}}\left(\sum_{i=1}^{n}\phi_{j}(z_{i})\epsilon_{i}\right)^{2}}\lVert (\hat{\beta}_{GL}-\beta)^{l} \rVert_{2}\\\label{eqgl}
&&+r\sum_{l=1}^{g}\omega_{l}\lVert \beta^{l} \rVert_{2}-r\sum_{l=1}^{g}\omega_{l}\lVert \hat{\beta}^{l}_{GL} \rVert_{2}.
\end{eqnarray}
Set  $Z_{l}=n^{-1}\sqrt{\sum_{j\in G_{l}}\left(\sum_{i=1}^{n}\phi_{j}(z_{i})\epsilon_{i}\right)^{2}},$ for $l\in \{1,\dots,g\}$ and the event
\begin{equation}\label{A}
\mathcal{A}=\bigcap_{l=1}^{g}\left\{Z_{l}\leq r\omega_{l}/2 \right\}.
\end{equation}
We state the result on event $\mathcal{A}$ and find an upper bound of $\mathbb{P}(\mathcal{A}^{c}).$\\
\textbf{On the event} $\mathcal{A}$ :
$$R(f_{\hat{\beta}_{GL}})-R(f_{0})\leq  R(f_{\beta})-R(f_{0})+ 
r\sum_{l=1}^{g}\omega_{l}\lVert (\hat{\beta}_{GL}-\beta)^{l} \rVert_{2}+
r\sum_{l=1}^{g}\omega_{l}\lVert \beta^{l} \rVert_{2}-r\sum_{l=1}^{g}\omega_{l}\lVert \hat{\beta}^{l}_{GL} \rVert_{2}.$$
This implies that
$$R(f_{\hat{\beta}_{GL}})-R(f_{0})\leq  R(f_{\beta})-R(f_{0})+ 2r\sum_{l=1}^{g}\omega_{l}\lVert \beta^{l} \rVert_{2}.$$
We conclude that on the event $\mathcal{A}$ we have
 $$R(f_{\hat{\beta}_{GL}})-R(f_{0})\leq \inf_{\beta\in \mathbb{R}^{p}}\left\{R(f_{\beta})-R(f_{0})+
2r\lVert \beta \rVert_{2,1}\underset{1\leq l \leq g}\max\omega_{l} \right\}.$$
We now come to the bound of $\mathbb{P}(A^{c})$ and write
 \begin{align}\label{formula}
  \mathbb{P}(\mathcal{A}^{c})&= \mathbb{P}\left(\overset{g}{\underset{l=1}{\bigcup}}\{\sqrt{\sum_{j\in G_{l}}\left(\sum_{i=1}^{n}\phi_{j}(z_{i})\epsilon_{i}\right)^{2}}> nr\omega_{l}/2\}\right)\\
 &\leq \sum_{l=1}^{g}\mathbb{P}\left(\sqrt{\sum_{j\in G_{l}}\left(\sum_{i=1}^{n}\phi_{j}(z_{i})\epsilon_{i}\right)^{2}}> nr\omega_{l}/2\right).
 \end{align}
 For $j\in G_{l}$ set $T_{j}^{l}=\sum_{i=1}^{n}\phi_{j}(z_{i})\epsilon_{i},$ we have  
 \begin{align*}
  \mathbb{P}(\mathcal{A}^{c})&\leq \sum_{l=1}^{g}\mathbb{P}\left(\sqrt{\sum_{j\in G_{l}}(T_{j}^{l})^{2}}> nr\omega_{l}/2\right)\\
  &= \sum_{l=1}^{g}\mathbb{P}\left(\sum_{j\in G_{l}}(T_{j}^{l})^{2}> (nr\omega_{l})^2/4\right).
\end{align*}
Using the fact that, for all $l\in\{1,\dots,g\}$
\begin{equation}
\left\{ \sum_{j\in G_{l}}(T_{j}^{l})^2> (nr\omega_{l})/4 \right\} \subset \underset{j\in G_{l}} \cup \left\{ (T_{j}^{l})^2> \frac{(nr\omega_{l})^2}{4|G_{l}|} \right\},
\end{equation}
it follows that
\begin{align*}
  \mathbb{P}(\mathcal{A}^{c})&\leq \sum_{l=1}^{g}\sum_{j\in G_{l}}\mathbb{P}\left(|T_{j}^{l}|> \frac{ nr\omega_{l}}{2\sqrt{|G_{l}|}}\right).
\end{align*}
For  $j\in G_{l},$  set $v_{j}^{l}=\sum_{i=1}^{n}\E(\phi_{j}^{2}\epsilon_{i}^{2}).$
 Since $\sum_{i=1}^{n}\phi_{j}^{2}(z_{i})\geqslant 4v_{j}^{l},$ we have
$$\mathbb{P}(|T_{j}^{l}|> \frac{nr\omega_{l}}{2\sqrt{|G_{l}|}})\leq 
\mathbb{P}\left(|T_{j}^{l}|> \sqrt{2v_{j}^{l}\left(x+\log{p}\right)}+\frac{c_{2}}{3}\left(x+\log{p}\right)\right),~~r\geq 1.$$
 By applying  Bernstein's inequality (see Lemma~\ref{Bern}) to the right hand side of the previous inequality we get
$$\mathbb{P}(|T_{j}^{l}|> \frac{n\omega_{l}}{2\sqrt{|G_{l}|}})\leq 2\exp\left(-x-\log{p}\right).$$
 It follows that
 \begin{align}\label{proA} 
 \mathbb{P}(\mathcal{A}^{c})&\leq \sum_{l=1}^{g}\sum_{j\in G_{l}}\mathbb{P}\left(|T_{j}^{l}|> \frac{n\omega_{l}}{2\sqrt{|G_{l}|}}\right)   \leq 2 \exp(-x).
 \end{align}
This ends the proof of the Theorem~\ref{theo02}. $\blacksquare$

\subsection{Proof of Theorem~\ref{theo21}}
 Fix an arbitrary $\beta\in\mathbb{R}^{p}$ such that $f_{\beta}\in \varGamma$. Set $\delta=W(\hat{\beta}_{GL}-\beta)$ where
$W=\mbox{Diag}(W_{1},\dots,W_{p})$ is a block diagonal matrix, with $W_{l}=\mbox{Diag}(\omega_{l},\dots,\omega_{l}).$
 Since
$\hat{\beta}_{GL}$ is the minimizer of $\hat{R}(f_{\beta})+r\sum_{l=1}^{g}\omega_{l}\lVert \beta^{l} \rVert_{2}$, we get
$$R(f_{\hat{\beta}_{GL}})-\frac{1}{n}\varepsilon^{T}X\hat{\beta}_{GL} +r\sum_{l=1}^{g}\omega_{l}\lVert \hat{\beta}^{l}_{GL} \rVert_{2}\leq R(f_{\beta})-\frac{1}{n}\varepsilon^{T}X\beta +
r\sum_{l=1}^{g}\omega_{l}\lVert \beta^{l} \rVert_{2}.$$
\textbf{On the event} $\mathcal{A}$ defined in~(\ref{A}), adding the term $\frac{r}{2}\sum_{l=1}^{g}\omega_{l}\lVert (\hat{\beta}_{GL}-\beta)^{l} \rVert_{2}$ to both sides of
 Inequality (\ref{eqgl}) yields to
\begin{align*}
 R(f_{\hat{\beta}_{GL}})+\frac{r}{2}\sum_{l=1}^{g}\omega_{l}\lVert (\hat{\beta}_{GL}-\beta)^{l} \rVert_{2}&\leq  R(f_{\beta})+
r\sum_{l=1}^{g}\omega_{l}(\lVert (\hat{\beta}_{GL}-\beta)^{l} \rVert_{2}-\lVert \hat{\beta}_{GL}^{l}\rVert_{2}+\lVert \beta^{l} \rVert_{2}).\\
\end{align*}
Since $ \lVert (\hat{\beta}_{GL}-\beta)^{l} \rVert_{2}-\lVert \hat{\beta}_{GL}^{l}\rVert_{2}+\lVert \beta^{l} \rVert_{2}=0$ for for $l \notin J(\beta)=J,$ we have
\begin{equation}\label{ll3}
 R(f_{\hat{\beta}_{GL}})-R(f_{0}) +\frac{r}{2}\sum_{l=1}^{g}\omega_{l}\lVert (\hat{\beta}_{GL}-\beta)^{l}\rVert_{2}\leq R(f_{\beta})-R(f_{0})+
2r\sum_{l\in J}\omega_{l}\lVert (\hat{\beta}_{GL}-\beta)^{l}\rVert_{2}.
\end{equation}
we get from Equation~(\ref{ll3}) 
that
\begin{equation}\label{deb11}
 R(f_{\hat{\beta}_{GL}})-R(f_{0})\leq R(f_{\beta})-R(f_{0})+2r\sum_{l\in J}\omega_{l}\lVert (\hat{\beta}_{GL}-\beta)^{l}\rVert_{2}
\end{equation}
 Consider separately  the two events :
\begin{equation*}\label{cas111}
 \mathcal{A}_{1}=\{2r\sum_{l\in J}\omega_{l}\lVert (\hat{\beta}_{GL}-\beta)^{l}\rVert_{2} \leq \eta(R(f_{\beta})-R(f_{0}))\},
\end{equation*}
and
\begin{equation}\label{cas211}
\mathcal{A}_{1}^{c}=\{ \eta(R(f_{\beta})-R(f_{0}))<2r\sum_{l\in J}\omega_{l}\lVert (\hat{\beta}_{GL}-\beta)^{l}\rVert_{2}\}.
\end{equation}
On the event $\mathcal{A}\cap \mathcal{A}_{1},$ we get from~(\ref{deb11})
\begin{equation}\label{eps11}
 R(f_{\hat{\beta}_{GL}})-R(f_{0})\leq (1+\eta)(R(f_{\beta})-R(f_{0})),
\end{equation}
and the result follows.
On the event $\mathcal{A}\cap \mathcal{A}_{1}^{c},$ all the following inequalities are valid. On one hand,
by applying Cauchy Schwarz inequality, we get from (\ref{deb11}) that
\begin{align}
 \notag R(f_{\hat{\beta}_{GL}})-R(f_{0})&\leq R(f_{\beta})-R(f_{0})+2r\sqrt{|J(\beta)|}\sqrt{\sum_{l\in J}\omega_{l}^{2}\lVert (\hat{\beta}_{GL}-\beta)^{l}\rVert_{2}^{2}}\\
&\leq R(f_{\beta})-R(f_{0})+2r\sqrt{|J(\beta)|}\lVert \delta_{J}\rVert_{2}. \label{eq41}
\end{align}
On the other hand we get from Equation~(\ref{ll3}) that
$$\frac{1}{2}\sum_{l=1}^{g}\omega_{l}\lVert (\hat{\beta}_{GL}-\beta)^{l}\rVert_{2}\leq R(f_{\beta})-R(f_{0})+2r\sum_{l\in J}\omega_{l}\lVert (\hat{\beta}_{GL}-\beta)^{l}\rVert_{2},$$
and using (\ref{cas211}) we obtain
$$\frac{1}{2}\sum_{l\in J}\omega_{l}\lVert (\hat{\beta}_{GL}-\beta)^{l}\rVert_{2}+ \frac{1}{2}\sum_{l\in J^{c}}\omega_{l}\lVert (\hat{\beta}_{GL}-\beta)^{l}\rVert_{2}
\leq \frac{2}{\eta} \sum_{l\in J}\omega_{l}\lVert (\hat{\beta}_{GL}-\beta)^{l}\rVert_{2}+2 \sum_{l\in J}\omega_{l}\lVert (\hat{\beta}_{GL}-\beta)^{l}\rVert_{2},$$
which implies $$\lVert \delta_{J^{c}} \rVert_{2,1}\leq (3+4/\eta)\lVert \delta_{J} \rVert_{2,1}.$$
 We can therefore apply Assumption~(\ref{AA2}) with  $a_{0}=3+4/\eta,$  and conclude that
\begin{equation}\label{eqmu11}
 \mu_{1}^{2} \lVert \delta_{J} \rVert_{2}^{2}\leq \frac{\lVert X\delta \rVert_{2}^{2}}{n}=\frac{1}{n}(\hat{\beta}_{GL}-\beta)^{T}WX^{T}XW(\hat{\beta}_{GL}-\beta)\leq 
(\underset{1\leq l \leq g}\max\omega_{l})^{2}\lVert f_{\hat{\beta}_{GL}}-f_{\beta}\rVert_{n}^{2}.
\end{equation}
Gathering Equations~(\ref{eq41}) and (\ref{eqmu11}) we get
\begin{align*}
R(f_{\hat{\beta}_{GL}})-R(f_{0})&\leq R(f_{\beta})-R(f_{0})+2r(\underset{1\leq l \leq g}\max\omega_{l})\sqrt{|J(\beta)|}\mu_{1}^{-1}\lVert f_{\hat{\beta}_{GL}}-f_{\beta}\rVert_{n}\\
&\leq R(f_{\beta})-R(f_{0})+2r(\underset{1\leq l \leq g}\max\omega_{l})\sqrt{|J(\beta)|}\mu_{1}^{-1}(\lVert f_{\hat{\beta}_{GL}}-f_{0}\rVert_{n}+\lVert f_{\beta}-f_{0}\rVert_{n}).
\end{align*}
We now use Lemma~\ref{l81} which compares excess risk to empirical norm.
\begin{lemma}\label{l81}
 Under assumptions~(\ref{A1}) and (\ref{A3}) we have
$$c_{0}\epsilon_{0}\lVert f_{\beta}-f_{0}\rVert_{n}^{2}\leq R(f_{\beta})-R(f_{0})\leq \frac{1}{4}c_{0}^{\prime}\lVert f_{\beta}-f_{0}\rVert_{n}^{2}.$$
where $c_{0}$ and $c_{0}^{\prime}$ are constants depending on $C_{0}$; and $\epsilon_{0}$ is a constant depending on $c_{1}$ and $c_{2}.$
\end{lemma}
(See the Appendix for the proof of Lemma~\ref{l81}).\\
Consequently 
\begin{eqnarray*}
R(f_{\hat{\beta}_{GL}})-R(f_{0})&\leq& R(f_{\beta})-R(f_{0})+\frac{2r(\underset{1\leq l \leq g}\max\omega_{l})\sqrt{|J(\beta)|}\mu_{1}^{-1}}{\sqrt{c_{0}\epsilon_{0}}}\sqrt{R(f_{\hat{\beta}_{GL}})-R(f_{0})}\\
& &+ \frac{2r(\underset{1\leq l \leq g}\max\omega_{l})\sqrt{|J(\beta)|}\mu_{1}^{-1}}{\sqrt{c_{0}\epsilon_{0}}}\sqrt{R(f_{\beta})-R(f_{0})}.
\end{eqnarray*}
Using inequality  $2uv<u^{2}/b+bv^{2}$ for all $b>1,$  with $u=r(\underset{1\leq l \leq g}\max\omega_{l})\frac{\sqrt{|J(\beta)|}\mu_{1}^{-1}}{\sqrt{c_{0}\epsilon_{0}}}$
and $v$ being either\\
  $\sqrt{R(f_{\hat{\beta}_{GL}})-R(f_{0})}$ or $\sqrt{R(f_{\beta})-R(f_{0})}$ we have
\begin{eqnarray*}
R(f_{\hat{\beta}_{GL}})-R(f_{0})&\leq& R(f_{\beta})-R(f_{0})+2b\left(\frac{r(\underset{1\leq l \leq g}\max\omega_{l})\sqrt{|J(\beta)|}\mu_{1}^{-1}}{\sqrt{c_{0}\epsilon_{0}}}\right)^{2}\\
& & +\frac{R(f_{\hat{\beta}_{GL}})-R(f_{0})}{b}+\frac{R(f_{\beta})-R(f_{0})}{b}.
\end{eqnarray*}
This implies that 
\begin{equation}\label{eqlem2}
R(f_{\hat{\beta}_{GL}})-R(f_{0})\leq \frac{b+1}{b-1}\left\{ R(f_{\beta})-R(f_{0})+\frac{2b^{2}r^{2}(\underset{1\leq l \leq g}\max\omega_{l})^{2}|J(\beta)|}{(b+1)\mu_{1}^{2} c_{0}\epsilon_{0}}\right\}.
\end{equation}
Now taking $b=1+2/\eta$  leads to 
\begin{equation}\label{eps12}
 R(f_{\hat{\beta}_{GL}})-R(f_{0})\leq (1+\eta)\left\{ R(f_{\beta})-R(f_{0})+\frac{c(\eta)r^{2}(\underset{1\leq l \leq g}\max\omega_{l})^{2}|J(\beta)|}{\mu_{1}^{2} c_{0}\epsilon_{0}}\right\}.
\end{equation}
According to Inequalities~(\ref{eps11}) and (\ref{eps12}) we conclude that on event $\mathcal{A},$ 
\begin{equation}\label{eps13}
 R(f_{\hat{\beta}_{GL}})-R(f_{0})\leq (1+\eta)\left\{ R(f_{\beta})-R(f_{0})+\frac{c(\eta)r^{2}(\underset{1\leq l \leq g}\max\omega_{l})^{2}|J(\beta)|}{\mu_{1}^{2} c_{0}\epsilon_{0}}\right\},
\end{equation}
where $c(\eta)=2(1+2/\eta)^{2}/(2+2/\eta).$ Inequality~(\ref{theo21eq1}) of the Theorem~\ref{theo21} follows.
Inequality~(\ref{theo2eq2}) follows from Lemma~\ref{l81}. This ends the proof of the Theorem~\ref{theo21} by considering (\ref{proA}).  $\blacksquare$
\subsection{Proof of Corollary~\ref{col2}}
Set $ \delta=W(\hat{\beta}_{GL}-\beta_{0}),$
 Line~(\ref{col21}) of Corollary~\ref{col2} follows directly from Equation~(\ref{eps13}) with $\beta=\beta_{0}$ and $\eta=1.$ Note that on the event $\mathcal{A}$ defined in~(\ref{A}), we have
\begin{equation}\label{k}
  \lVert \delta_{J(\beta_{0})^{c}}\rVert_{2,1}\leq 3\lVert \delta_{J(\beta_{0})}\rVert_{2,1}.
\end{equation}
Indeed, since $\hat{\beta}_{GL}$ is the
minimizer of $\hat{R}(f_{\beta})+r\sum_{l=1}^{g}\omega_{l}\lVert \beta^{l} \rVert_{2},$
$$R(f_{\hat{\beta}_{GL}})-R(f_{\beta_{0}})+r\sum_{l=1}^{g}\omega_{l}\lVert \hat{\beta}_{GL}^{l} \rVert_{2}\leq \frac{1}{n}\varepsilon^{T}X(\hat{\beta}_{GL}-\beta_{0})+r\sum_{l=1}^{g}\omega_{l}\Vert \beta^{l}_{0}\rVert_{2}$$
which implies
$$r\Vert W\hat{\beta}_{GL}\rVert_{2,1}\leq 
\sum_{l=1}^{g}\frac{1}{n}\sqrt{\sum_{j\in G_{l}}\left(\sum_{i=1}^{n}(z_{ij})\epsilon_{i}\right)^{2}}\lVert (\hat{\beta}_{GL}-\beta_{0})^{l} \rVert_{2}
+r\Vert W\beta_{0}\rVert_{2,1}$$
On the event $A$ we have
\begin{eqnarray*}
  \Vert W(\hat{\beta}_{GL})_{J(\beta_{0})}\rVert_{2,1}+\Vert W(\hat{\beta}_{GL})_{J^{c}(\beta_{0})}\rVert_{2,1}&\leq
&\frac{1}{2}(\Vert W(\hat{\beta}_{GL}-\beta_{0})_{J(\beta_{0})}\rVert_{2,1}+\Vert W(\hat{\beta}_{GL})_{J^{c}(\beta_{0})}\rVert_{2,1})\\
& &+\Vert W(\beta_{0})_{J(\beta_{0})}\rVert_{2,1}. 
\end{eqnarray*}
This yields to~(\ref{k}).
Line~(\ref{col22}) follows from Line~(\ref{col21}) by applying Lemma~\ref{l81}. Line~(\ref{col23}) follows from Line~(\ref{col22})  by using Equation~(\ref{eqmu11}) 
and $\lVert \delta\rVert_{2,1}^{2}\leq 16 s \lVert \delta_{J(\beta_{0})}\rVert_{2}^{2} $.
Line~(\ref{col24}) is the consequence of the Lemma~\ref{l9}
with $a_{l}=\lVert(\hat{\beta}_{GL}-\beta_{0})^{l}\rVert_{2}$ and 
$$b_{1}=\frac{12rs\left(\underset{1\leq l \leq g}\max\omega_{l}\right)^{2} }{\mu^{2}(s,3) c_{0}\epsilon_{0}(\underset{1\leq l \leq g}\min\omega_{l})}.~~~  \blacksquare$$
\subsection{Proof of Theorem~\ref{theo22}}
On the event $\mathcal{A}$ defined in~(\ref{A}), using Inequality~(\ref{eqgl}) with $\beta=\beta_0$ yields
\begin{equation}\label{eqq0}
R(f_{\hat{\beta}_{GL}})-R(f_{\beta_{0}})\leq  
\sum_{l=1}^{g}\frac{3r\omega_{l}}{2}\lVert (\hat{\beta}_{GL}-\beta_{0})^{l} \rVert_{2}.
\end{equation}
By Lemma~\ref{l1} we have,
\begin{equation}\label{eq0}
 \frac{\langle h,h \rangle_{f_{\beta_{0}}}}{\lVert h\rVert_{\infty}^{2}}(\exp(-\lVert h\rVert_{\infty})+\lVert h\rVert_{\infty}-1)\leq R(f_{\hat{\beta}_{GL}})-R(f_{\beta_{0}})
\end{equation}
where $$ h(z_{i})=(f_{\hat{\beta}_{GL}}-f_{\beta_{0}})(z_{i})=\sum_{l=1}^{g}\sum_{j\in G_{l}}(\hat{\beta}_{GL,j}-\beta_{0 j})z_{ij}.$$ 
One can easily verify that$\lVert h\rVert_{\infty}\leq v\lVert \delta^{\prime} \rVert_{2,1}$ with $\delta^{\prime}=\hat{\beta}_{GL}-\beta_{0}.$
Equation~(\ref{eq0}) and the decreasing of $t\mapsto \frac{\exp(-t)+t-1}{t^{2}}$ lead to
$$\frac{\delta^{'T}X^{T} DX\delta^{\prime}}{ n(v\lVert \delta^{\prime} \rVert_{2,1})^{2}}
(\exp(-v \lVert \delta^{\prime} \rVert_{2,1})+v\lVert \delta^{\prime} \rVert_{2,1}-1)\leq R(f_{\hat{\beta}_{GL}})-R(f_{\beta_{0}}).$$
Now, Inequality~(\ref{k}) implies
$$\lVert \delta^{\prime}_{J(\beta_{0})^{c}}\rVert_{2,1}\leq 3\frac{\left(\underset{1\leq l \leq g}\max\omega_{l}\right)}{\underset{1\leq l \leq g}\min\omega_{l}}\lVert \delta^{\prime}_{J(\beta_{0})}\rVert_{2,1}.$$
 We can therefore apply Assumption~(\ref{AAA2}) with  $a_{0}=3(\underset{1\leq l \leq g}\max\omega_{l})/\underset{1\leq l \leq g}\min\omega_{l}$ and get that
$$\frac{\mu_{2}^{2}\lVert\delta^{\prime}_{J} \rVert_{2}^{2}}{ v^{2}\lVert \delta^{\prime} \rVert_{2,1}^{2}}
(\exp(-v \lVert \delta^{\prime} \rVert_{2,1})+v\lVert \delta^{\prime} \rVert_{2,1}-1)\leq R(f_{\hat{\beta}_{GL}})-R(f_{\beta_{0}}).$$
We can  use that $\lVert \delta^{\prime}\rVert_{2,1}^{2}\leq (1+a_{0})^{2} |J| \lVert\delta^{\prime}_{J}\rVert_{2}^{2},$ with $J=J(\beta_{0})$ to write
$$\frac{\mu_{2}^{2}}{(1+a_{0})^{2} |J| v^{2}}
(\exp(-v \lVert \delta^{\prime} \rVert_{2,1})+v\lVert \delta^{\prime} \rVert_{2,1}-1)\leq R(f_{\hat{\beta}_{GL}})-R(f_{\beta_{0}}).$$
According to Equation~(\ref{eqq0}) we have
\begin{equation}\label{eqq1}
 \exp(-v \lVert \delta^{\prime} \rVert_{2,1})+v\lVert \delta^{\prime} \rVert_{2,1}-1\leq \frac{3r(1+a_{0})^{2} \left(\underset{1\leq l \leq g}\max\omega_{l}\right)v^{2}|J|}{2\mu_{2}^{2}}\Vert \delta^{\prime}\rVert_{2,1}.
\end{equation}
Now, a short calculation shows that for all $ a\in (0,1],$
\begin{equation}
 e^{\frac{-2a}{1-a}}+(1-a)\frac{2a}{1-a}-1\geqslant 0
\end{equation}
Set $a=v\lVert \delta^{\prime} \rVert_{2,1}/(v\lVert \delta^{\prime} \rVert_{2,1}+2).$ Thus $v\lVert \delta^{\prime} \rVert_{2,1}=2a/(1-a)$
and we have
\begin{equation}\label{kk}
  e^{-v\lVert \delta^{\prime} \rVert_{2,1}}+v\lVert \delta^{\prime} \rVert_{2,1}-1  \geqslant \frac{v^{2}\lVert \delta^{\prime} \rVert_{2,1}^{2}}{v\lVert \delta^{\prime} \rVert_{2,1}+2}.
\end{equation}
This implies using Equation~(\ref{eqq1}) that 
$$v\lVert \delta^{\prime} \rVert_{2,1}\leq \frac{3r(1+a_{0})^{2} \left(\underset{1\leq l \leq g}\max\omega_{l}\right)|J| v/\mu_{2}^{2}}{1-3r(1+a_0)^{2} \left(\underset{1\leq l \leq g}\max\omega_{l}\right)|J| v/2\mu_{2}^{2}}.$$
Now if  $r(1+a_0)^{2}\underset{1\leq l \leq g}\max\omega_{l}\leq \frac{\mu_{2}^{2}}{3v|J|},$ we have $v\lVert \delta^{\prime} \rVert_{2,1}\leq 2$
and consequently 
$$\frac{\exp(-v \lVert \delta^{\prime} \rVert_{2,1})+v\lVert \delta^{\prime} \rVert_{2,1}-1}{v^{2}\lVert \delta^{\prime} \rVert_{2,1}^{2}}\geqslant 1/4.$$
Now,  Inequality~(\ref{eqq1}) implies
$$\lVert \delta^{\prime} \rVert_{2,1}\leq \frac{6(1+a_0)^{2} |J|r\left(\underset{1\leq l \leq g}\max\omega_{l}\right)}{\mu_{2}^{2}}.$$
This proves the Line~(\ref{theo222}). Line~(\ref{theo221}) follows from (\ref{theo222}) by using Inequality~(\ref{eqq0}).
Line~(\ref{theo223}) is the consequence of Lemma~\ref{l9}
taking $a_{l}=\lVert(\hat{\beta}_{GL}-\beta_{0})^{l}\rVert_{2}$ and $b_{1}=6(1+a_0)^{2}|J|r(\underset{1\leq l \leq g}\min\omega_{l})/\mu_{2}^{2}(s,3)$.
Line~(\ref{theo224}) follows from Line~(\ref{theo221}) and Inequality~(\ref{eq0}). $\blacksquare$
\subsection{Proof of Theorem~\ref{theo0}}
Note that Lasso can be derived by Group Lasso by taking one predictor per group i.e $p=g$ and $G_j=\{j\}$ for $j\in\{1,\dots,p\}.$
This implies, using (\ref{eqgl}) that
$$R(f_{\hat{\beta}_{L}})-R(f_{0}) \leq
 R(f_{\beta})-R(f_{0})+ \sum_{j=1}^{p}\left|\frac{1}{n}\sum_{i=1}^{n}\phi_{j}(z_{i})\varepsilon_{i}\right||\hat{\beta}_{L,j}-\beta_{j}| +r\sum_{j=1}^{p}\omega_{j}|\beta_{j}|-r\sum_{j=1}^{p}\omega_{j}|\hat{\beta}_{L,j}|.$$
For $1\leq j\leq p,$ set $S_{j}=\sum_{i=1}^{n}\phi_{j}(z_{i})\varepsilon_{i}$ and let us denote by $E,$ the event
\begin{equation}\label{E}
 E=\bigcap_{j=1}^{p}\left\{|S_{j}|\leq nr\omega_{j}/2 \right\}.
\end{equation}
We state the results on the event $E$ and then find an upper bound of  $\mathbb{P}(E^{c})$.\\
\textbf{On the event} $E$ :
\begin{align*}
R(f_{\hat{\beta}_{L}})-R(f_{0}) &\leq R(f_{\beta})-R(f_{0})+
r\sum_{j=1}^{p}\omega_{j}|\hat{\beta}_{L,j}-\beta_{j}| +r\sum_{j=1}^{p}\omega_{j}|\beta_{j}|-r\sum_{j=1}^{p}\omega_{j}|\hat{\beta}_{L,j}|\\
 &\leq R(f_{\beta})-R(f_{0})+2r\sum_{j=1}^{p}\omega_{j}|\beta_{j}|.
\end{align*}
We conclude that on the event $E$ we have
\begin{equation*}
 R(f_{\hat{\beta}_{L}})-R(f_{0})\leq \inf_{\beta\in \mathbb{R}^{p}}\left\{R(f_{\beta})-R(f_{0})+
   2r\lVert \beta \rVert_{1}\underset{1\leq j \leq p}\max\omega_{j}\right\}.
\end{equation*}
Now we are going to find  an upper bound of  $\mathbb{P}(E^{c}) :$
\begin{align*} 
\mathbb{P}(E^{c})&\leq\mathbb{P}\left(\overset{p}{\underset{j=1}{\bigcup}}\{|\sum_{i=1}^{n}\phi_{j}(z_{i})(Y_{i}-\E(Y_{i}))|> r\omega_{j} n/2\}\right)\\
&\leq \sum_{j=1}^{p}\mathbb{P}(|\sum_{i=1}^{n}\phi_{j}(z_{i})(Y_{i}-\E(Y_{i}))|>r\omega_{j} n/2).
\end{align*}
 For  $j\in\{1,\dots,p\},$  set $v_{j}=\sum_{i=1}^{n}\E(\phi_{j}^{2}\epsilon_{i}^{2}).$
Since $\sum_{i=1}^{n}\phi_{j}^{2}(z_{i})\geqslant 4 v_{j},$ we have
$$\mathbb{P}(|S_{j}|> nr\omega_{j}/2)\leq \mathbb{P}\left(|S_{j}|> \sqrt{2v_{j}(x+\log{p})}+\frac{c_{2}}{3}(x+\log{p})\right),~~r\geq 1.$$
By applying  Bernstein's inequality (see ~\cite{boucheron, massart2007}) to the right hand side of the previous inequality we get
$$\mathbb{P}(|S_{j}|> nr\omega_{j}/2)\leq 2\exp(-x-\log{p}).$$
It follows that
\begin{align}\label{proE1} 
\mathbb{P}(E^{c})&\leq \sum_{j=1}^{p}\mathbb{P}(|S_{j}|>r\omega_{j} n/2)\leq 2\exp(-x).
\end{align}
When  $\omega_{j}=1,~\mbox{for ~all} ~j \in \{1,\dots,p\}$ and $r=A\sqrt{\frac{\log{p}}{n}}$, we apply Hoeffding's inequality (see ~\cite{boucheron, massart2007}). This leads to
\begin{align} 
\notag \mathbb{P}(E^{c})&=\mathbb{P}\left(\overset{p}{\underset{j=1}{\bigcup}}\{|\sum_{i=1}^{n}\phi_{j}(z_{i})(Y_{i}-\E(Y_{i}))|> r n/2\}\right)\\
\notag&\leq \sum_{j=1}^{p}\mathbb{P}(|\sum_{i=1}^{n}\phi_{j}(z_{i})(Y_{i}-\E(Y_{i}))|>r n/2)\\
 \label{proE2}&\leq 2p\exp\left(-\frac{2(r n/2)^{2}}{\sum_{i=1}^{n}2c_{2}}\right)=2p\exp\left(-\frac{r^{2}n}{4c_{2}}\right)=2p^{1-\frac{A^{2}}{4c_{2}}}.
\end{align}
This ends the proof of Theorem~\ref{theo0}. $\blacksquare$
\subsection{Proof of Theorem~\ref{theo1}}
 Fix an arbitrary $\beta\in\mathbb{R}^{p}$ such that $f_{\beta}\in \varGamma,$ and set $\delta=W(\hat{\beta}_{L}-\beta)$, where $W=\mbox{Diag}(w_{1},\dots,w_{p}).$ 
It follows from Inequality~(\ref{eps13}) that
\begin{equation}\label{eps3}
 R(f_{\hat{\beta}_{L}})-R(f_{0})\leq (1+\eta)\left\{ R(f_{\beta})-R(f_{0})+\frac{c(\eta)r^{2}\left(\underset{1\leq j \leq p}\max\omega_{j}\right)^{2}|K(\beta)|}{\mu^{2} c_{0}\epsilon_{0}}\right\},
\end{equation}
where $c(\eta)=2(1+2/\eta)^{2}\Large{/}(2+2/\eta).$
This ends the proof of Inequality~(\ref{theo1eq1}) of the Theorem~\ref{theo1}. Inequality~(\ref{theo1eq2}) follows from Lemma~\ref{l81}.
To prove Inequalities~(\ref{theo1eq3}) and (\ref{theo1eq4}) we just replace $\omega_{j}$ by  $A\sqrt{\frac{\log{p}}{n}}.$\\
This ends the proof of the Theorem~\ref{theo1} by using (\ref{proE1}) and (\ref{proE2}). $\blacksquare$
\subsection{Proof of Corollary~\ref{col1}}
Set $ \delta=W(\hat{\beta}_{L}-\beta_{0}).$
The result~(\ref{col1eq1}) directly comes by taking $\beta=\beta_{0}$ and $\eta=2$ in (\ref{eps3}).
Note that, on the event $E$ defined in~(\ref{E}), we have
\begin{equation}\label{J}
\lVert \delta_{K(\beta_{0})^{c}}\rVert_{1}\leq 3\lVert \delta_{K(\beta_{0})}\rVert_{1}.
\end{equation}
Indeed, since $\hat{\beta}_{L}$ is the
minimizer of $\hat{R}(f_{\beta})+r\sum_{j=1}^{p}\omega_{j}|\beta_{j}|,$ then
$$R(f_{\hat{\beta}_{L}})-R(f_{\beta_{0}})+r\sum_{j=1}^{p}\omega_{j}|\hat{\beta}_{L,j}|\leq 
\frac{1}{n}\varepsilon^{T}X(\hat{\beta}_{L}-\beta_{0}) +r\sum_{j=1}^{p}\omega_{j}|\beta_{0j}|,$$
which implies that
$$r\Vert W\hat{\beta}_{L}\rVert_{1}\leq 
\sum_{j=1}^{p}\left|\frac{1}{n}\sum_{i=1}^{n}\phi_{j}(z_{i})\varepsilon_{i}\right||\hat{\beta}_{L,j}-\beta_{j}| 
+r\Vert W\beta_{0}\rVert_{1}.$$
On the event $E$ we have
\begin{eqnarray*}
  \Vert W(\hat{\beta}_{L})_{K(\beta_{0})}\rVert_{1}+\Vert W(\hat{\beta}_{L})_{K^{c}(\beta_{0})}\rVert_{1}&\leq&\frac{1}{2}(\Vert W(\hat{\beta}_{L}-\beta_{0})_{K(\beta_{0})}\rVert_{1}+\Vert W(\hat{\beta}_{L})_{K^{c}(\beta_{0})}\rVert_{1})\\
& &+\Vert W(\beta_{0})_{K(\beta_{0})}\rVert_{1}. 
\end{eqnarray*}
Thus (\ref{J}) follows.
Line~(\ref{col1eq2}) follows from Line~(\ref{col1eq1}) by applying Lemma~\ref{l81}. Line~(\ref{col1eq3})  follows from Line(\ref{col1eq2}) by using Inequality~(\ref{eqmu11})
and $\lVert \delta\rVert_{1}^{2}\leq 16 s \lVert \delta_{K(\beta_{0})}\rVert_{2}^{2} $.
The last line follows from Lemma~\ref{l9} in Appendix
with $a_{j}=|\hat{\beta}_{L, j}-\beta_{0j}|$ and $$b_{1}= \frac{12 sr \left(\underset{1\leq j \leq p}\max\omega_{j}\right)^{2}}{\mu^{2}(s,3) c_{0}\epsilon_{0}\left(\underset{1\leq j \leq p}\min\omega_{j}\right)}. ~~~\blacksquare$$

\section{Appendix}
The proof of Lemma~\ref{l81} are based on property of self concordant function (see for instance \cite{nesterov1987}), \textit{i.e.}, the functions whose third
derivatives are controlled by their second derivatives. A one-dimensional, convex function $g$ is called self concordant if
$$|g^{'''}(x)|\leq C g^{''}(x)^{3/2}.$$
The function we use $(g(t)=\hat{R}(g+th))$ is not really  self concordant but we can bound his third derivative by the second derivative times a constant.
Our results on self-concordant functions are based on the ones of \cite{bach10}.  He has used and extended tools from convex optimization and self-concordance to provide simple
extensions of theoretical results for the square loss to logistic loss. We use the same kind of arguments and state some relations between 
excess risk and prediction loss in the context of nonparametric logistic model, where $f_{0}$ is not necessarily  linear as assumed in \cite{bach10}.
Precisely we extend Proposition 1 in \cite{bach10} to the functions which are not necessarily  linear (see Lemma~\ref{l1}). This allows us to establish Lemma~\ref{l81}. 
\begin{lemma}\label{l1}
 For all $h,f : \mathbb{R}^{d}\rightarrow \mathbb{R}$, we have 
\begin{equation}
\frac{\langle h,h \rangle_{f}}{\lVert h\rVert_{\infty}^{2}}(\exp(-\lVert h\rVert_{\infty})+\lVert h\rVert_{\infty}-1)\leq R(f+h)-R(f)+(q_{f}-q_{f_{0}})(h),
\end{equation}
\begin{equation}
 R(f+h)-R(f)+(q_{f}-q_{f_{0}})(h)\leq \frac{\langle h,h \rangle_{f}}{\lVert h\rVert_{\infty}^{2}}(\exp(\lVert h\rVert_{\infty})-\lVert h\rVert_{\infty}-1),
\end{equation}
and
\begin{equation}
 \langle h,h\rangle_{f}e^{-\lVert h\rVert_{\infty}} \leq \langle h,h\rangle_{f+h}\leq \langle h,h\rangle_{f} e^{\lVert h\rVert_{\infty}}.
\end{equation}

\end{lemma}
\section{Proof of Lemma~\ref{l1}}
 We use the following lemma (see~\cite{bach10} Lemma 1) that we recall here :
\begin{lemma}\label{bach}
 Let $g$ be a convex three times differentiable function $g :\mathbb{R} \rightarrow \mathbb{R} $ such that for all $t\in \mathbb{R}$ $|g^{'''}(t)|\leq Sg^{''}(t)$, 
for some $S\geq 0$. Then , for all $t\geq0$ :
\begin{equation}
 \frac{g^{''}(0)}{S^{2}}(\exp(-St)+St-1)\leq g(t)-g(0)-g^{\prime}(0)t\leq \frac{g^{''}(0)}{S^{2}}(\exp(St)-St-1) .
\end{equation}
\end{lemma}
We refer to Appendix A of \cite{bach10} for the proof of this lemma. \\
Set $$g(t)=\hat{R}(f+th)=\frac{1}{n}\sum_{i=1}^{n}l((f+th)(z_{i}))-Y_{i}(f+th)(z_{i}),~~~f,h \in H,$$ 
where $l(u)=\log(1+\exp(u)).$ A short calculation leads to $l^{\prime}(u)=\pi(u)$, $l^{''}(u)=\pi(u)(1-\pi(u))$, $l^{'''}(u)=\pi(u)[1-\pi(u)][1-2\pi(u)].$
It follows that  
$$ g^{''}(t)=\frac{1}{n}\sum_{i=1}^{n}h^{2}(z_{i})l^{''}((f+th)(z_{i}))=\langle h,h\rangle_{f+th},$$ and 
$$g^{'''}(t)=\frac{1}{n}\sum_{i=1}^{n}h^{3}(z_{i})l^{'''}((f+th)(z_{i})).$$ 
Since
$l^{'''}(u)\leq l^{''}(u)$ we have,
\begin{align*}
| g^{'''}(t)|& =\left|\frac{1}{n}\sum_{i=1}^{n}h^{3}(z_{i})l^{'''}((f+th)(z_{i}))\right|\\
& \leq\frac{1}{n}\sum_{i=1}^{n}h^{2}(z_{i})l^{''}((f+th)(z_{i}))\lVert h\rVert_{\infty}=\lVert h\rVert_{\infty}g^{''}(t).
 \end{align*}
We now apply Lemma~\ref{bach} to $g(t)$ with $S=\lVert h\rVert_{\infty}$, taking $t=1.$  Using Equation~(\ref{equal}) 
we get the  first and second inequality of Lemma~\ref{l1}.
Now by considering $g(t)=\langle h,h\rangle_{f+th}$, a short calculation leads to $|g'(t)|\leq \lVert h\rVert_{\infty}g(t)$ which implies
$g(0)e^{-\lVert h\rVert_{\infty}t}\leq g(t)\leq g(0)e^{\lVert h\rVert_{\infty}t}.$ 
 By applying the last inequality to $g(t)$, and taking $t=1$ we get the third inequality of Lemma~\ref{l1}.

\subsection{Proof of Lemma~\ref{l81} }
Set $h_{0}=f_{\beta}-f_{0}$ from Lemma~\ref{l1} below, 
$$\frac{\langle h_{0},h_{0} \rangle_{f_{0}}}{\lVert h_{0}\rVert_{\infty}^{2}}(\exp(-\lVert h_{0}\rVert_{\infty})+
\lVert h_{0}\rVert_{\infty}-1)\leq R(f_{\beta})-R(f_{0}).$$
Using Assumptions~(\ref{A3}), ~(\ref{A1}) and the decreasing of $t\mapsto \frac{\exp(-t)+t-1}{t^{2}},$ we claim that there exists $c_{0}=c_{0}(C_{0},c_{1})>0$ such that
$$c_{0}\leq \frac{\exp(-\lVert h_{0}\rVert_{\infty})+\lVert h_{0}\rVert_{\infty}-1)}{\lVert h_{0}\rVert_{\infty}^{2}}.$$ According to Assumption~(\ref{A1}), 
there exists $0\leq \epsilon_{0}\leq 1/2$ such that for $1\leq i\leq n$
 $$\epsilon_{0}\leq \pi(f_{0}(z_{i}))(1-\pi(f_{0}(z_{i})))\leq 1-\epsilon_{0}.$$
The proof of the left hand side of Lemma~\ref{l81} follows from the fact that
$\epsilon_{0}\lVert h_{0}\rVert_{n}^{2}\leq \langle h_{0},h_{0} \rangle_{f_{0}}.$ 
From the second line of Lemma~\ref{l1} we have 
$$R(f_{\beta})-R(f_{0})\leq \frac{\langle h_{0},h_{0} \rangle_{f_{0}}}{\lVert h_{0}\rVert_{\infty}^{2}}(\exp(\lVert h_{0}\rVert_{\infty})-\lVert h_{0}\rVert_{\infty}-1).$$
Using assumption~(\ref{A3}) and increasing  of $t\mapsto \frac{\exp(t)-t-1}{t^{2}}$ thus there exists $c_{0}^{\prime}=c_{0}^{\prime}(C_{0},c_{1})>0$ such that
\begin{align*}
 R(f_{\beta})-R(f_{0})&\leq c_{0}^{\prime} \langle h_{0},h_{0} \rangle_{f_{0}}\\
&\leq c_{0}^{\prime}\frac{1}{4}\lVert h_{0}\rVert_{n}^{2}.
\end{align*}
This end the proof of the right hand side of the Lemma~\ref{l81}.
\begin{lemma}\label{l9}
 If we assume that $\sum_{i=1}^{p}a_{j}\leq b_{1}$ with $a_{j}>0,$ this implies that $\sum_{i=1}^{p}a^{q}_{j}\leq b_{1}^{q},$ with $1\leq q \leq 2$.
\end{lemma}
\subsection{Proof of Lemma~\ref{l9}}
We start by writing 
\begin{align*}
 \sum_{i=1}^{p}a^{q}_{j}&=\sum_{i=1}^{p}a^{2-q}_{j}a^{2q-2}_{j}\\
&\leq \left(\sum_{i=1}^{p}a_{j}\right)^{2-q}\left(\sum_{i=1}^{p}a^{2}_{j} \right)^{q-1}.
\end{align*}
Since $\sum_{i=1}^{p}a^{2}_{j}\leq \left(\sum_{i=1}^{p}a_{j}\right)^{2}\leq b_{1}^{2}$, thus
\begin{equation}\label{in1}
 \sum_{i=1}^{p}a^{q}_{j}\leq b_{1}^{2-q}b_{1}^{2q-2}=b_{1}^{q}.
\end{equation}
This ends the proof.
\begin{lemma}[Bernstein's inequality]\label{Bern}
 Let $X_{1},\dots,X_{n}$ be independent real valued random variables such that for all $i\leq n,$ $X_{i}\leq b$ almost surely, then 
we have
$$\mathbb{P}\left[\left|\sum_{i=1}^{n}X_{i}-\mathbb{E}(X_{i})\right|\geqslant \sqrt{2vx}+bx/3\right]\leq 2\exp(-x),$$
where $v=\sum_{i=1}^{n}\mathbb{E}(X_{i}^{2}).$
\end{lemma}
This lemma is obtain by gathering  Proposition 2.9  and  inequality (2.23) from ~\cite{massart2007}. 

\begin{lemma}[Hoeffding's inequality]\label{hoef}
 Let $X_{1},\dots,X_{n}$ be independent random variables such that $X_{i}$ takes its values in $[a_{i}, b_{i}]$ almost surely for all $i\leq n.$
Then for any positive x, we have  
$$\mathbb{P}\left[\left|\sum_{i=1}^{n}X_{i}-\mathbb{E}(X_{i})\right|\geqslant x\right]\leq 2\exp(-\frac{2x^{2}}{\sum_{i=1}^{n}(b_{i}-a_{i})^{2}}).$$
\end{lemma}
This lemma is a consequence of Proposition 2.7 in~\cite{massart2007}.
  \bibliographystyle{plain}
 \bibliography{biblioM_JMVA}

\def\cprime{$'$}
\begin{thebibliography}{10}

\bibitem{akaike}
H.~Akaike.
\newblock Information theory and an extension of the maximum likelihood
  principle.
\newblock In {\em Second {I}nternational {S}ymposium on {I}nformation {T}heory
  ({T}sahkadsor, 1971)}, pages 267--281. Akad\'emiai Kiad\'o, Budapest, 1973.

\bibitem{bach10}
Francis Bach.
\newblock Self-concordant analysis for logistic regression.
\newblock {\em Electronic Journal of Statistics}, 4:384--414, 2010.

\bibitem{bartlett2012}
Peter~L Bartlett, Shahar Mendelson, and Joseph Neeman.
\newblock L1-regularized linear regression: persistence and oracle
  inequalities.
\newblock {\em Probability theory and related fields}, 154(1-2):193--224, 2012.

\bibitem{BRT}
Peter~J. Bickel, Ya'acov Ritov, and Alexandre~B. Tsybakov.
\newblock Simultaneous analysis of lasso and {D}antzig selector.
\newblock {\em Annals of Statistics}, 37(4):1705--1732, 2009.

\bibitem{boucheron}
S.~Boucheron, G.~Lugosi, and O.~Bousquet.
\newblock Concentration inequalities.
\newblock {\em Advanced Lectures on Machine Learning}, pages 208--240, 2004.

\bibitem{bunea2007sparsity}
Florentina Bunea, Alexandre Tsybakov, and Marten Wegkamp.
\newblock Sparsity oracle inequalities for the {L}asso.
\newblock {\em Electronic Journal of Statistics}, 1:169--194, 2007.

\bibitem{bunea2006aggregation}
Florentina Bunea, Alexandre~B. Tsybakov, and Marten~H. Wegkamp.
\newblock Aggregation and sparsity via {$l_1$} penalized least squares.
\newblock In {\em Learning theory}, volume 4005 of {\em Lecture Notes in
  Comput. Sci.}, pages 379--391. Springer, Berlin, 2006.

\bibitem{bunea2007aggregation}
Florentina Bunea, Alexandre~B. Tsybakov, and Marten~H. Wegkamp.
\newblock Aggregation for {G}aussian regression.
\newblock {\em Annals of Statistics}, 35(4):1674--1697, 2007.

\bibitem{chesneau}
C.~Chesneau and M.~Hebiri.
\newblock Some theoretical results on the grouped variables lasso.
\newblock {\em Mathematical Methods of Statistics}, 17(4):317--326, 2008.

\bibitem{Fried}
J.~Friedman, T.~Hastie, and R.~Tibshirani.
\newblock Regularization paths for generalized linear models via coordinate
  descent.
\newblock {\em Journal of statistical software}, 33(1):1, 2010.

\bibitem{manuel}
Manuel Garcia-Magari{\~n}os, Anestis Antoniadis, Ricardo Cao, and Wenceslao
  Gonz{\'a}lez-Manteiga.
\newblock Lasso logistic regression, {GS}oft and the cyclic coordinate descent
  algorithm: application to gene expression data.
\newblock {\em Stat. Appl. Genet. Mol. Biol.}, 9:Art. 30, 30, 2010.

\bibitem{non_par_logist}
T~Hastie.
\newblock Non-parametric logistic regression.
\newblock {\em SLAC PUB-3160, June}, 1983.

\bibitem{huang2010}
J.~Huang, J.L. Horowitz, and F.~Wei.
\newblock Variable selection in nonparametric additive models.
\newblock {\em Annals of statistics}, 38(4):2282, 2010.

\bibitem{huang2008}
J.~Huang, S.~Ma, and CH~Zhang.
\newblock The iterated lasso for high--dimensional logistic regression.
\newblock {\em Technical Report 392}, 2008.

\bibitem{james2009dasso}
Gareth~M James, Peter Radchenko, and Jinchi Lv.
\newblock Dasso: connections between the dantzig selector and lasso.
\newblock {\em Journal of the Royal Statistical Society: Series B (Statistical
  Methodology)}, 71(1):127--142, 2009.

\bibitem{KF}
Keith Knight and Wenjiang Fu.
\newblock Asymptotics for lasso-type estimators.
\newblock {\em Annals of Statistics}, 28(5):1356--1378, 2000.

\bibitem{lounici2009}
K.~Lounici, M.~Pontil, A.B. Tsybakov, and S.~Van De~Geer.
\newblock Taking advantage of sparsity in multi-task learning.
\newblock {\em In COLT'09}, 2009.

\bibitem{lounici2011}
Karim Lounici, Massimiliano Pontil, Sara van~de Geer, and Alexandre~B.
  Tsybakov.
\newblock Oracle inequalities and optimal inference under group sparsity.
\newblock {\em Annals of Statistics}, 39(4):2164--2204, 2011.

\bibitem{massart2007}
Pascal Massart.
\newblock {\em Concentration inequalities and model selection}, volume 1896 of
  {\em Lecture Notes in Mathematics}.
\newblock Springer, Berlin, 2007.
\newblock Lectures from the 33rd Summer School on Probability Theory held in
  Saint-Flour, July 6--23, 2003, With a foreword by Jean Picard.

\bibitem{masmen}
Pascal Massart and Caroline Meynet.
\newblock The {L}asso as an {$\ell_1$}-ball model selection procedure.
\newblock {\em Electronic Journal of Statistics}, 5:669--687, 2011.

\bibitem{mcauley2005}
J.~McAuley, J.~Ming, D.~Stewart, and P.~Hanna.
\newblock Subband correlation and robust speech recognition.
\newblock {\em Speech and Audio Processing, IEEE Transactions on},
  13(5):956--964, 2005.

\bibitem{meier2009}
L.~Meier, S.~Van De~Geer, and P.~B{\"u}hlmann.
\newblock High-dimensional additive modeling.
\newblock {\em Annals of Statistics}, 37(6B):3779--3821, 2009.

\bibitem{mal}
Lukas Meier, Sara van~de Geer, and Peter B{\"u}hlmann.
\newblock The group {L}asso for logistic regression.
\newblock {\em Journal of the Royal Statistical Society Series B},
  70(1):53--71, 2008.

\bibitem{meinshausen2006}
Nicolai Meinshausen and Peter B{\"u}hlmann.
\newblock High-dimensional graphs and variable selection with the lasso.
\newblock {\em Annals of Statistics}, 34(3):1436--1462, 2006.

\bibitem{meinshausen2009lasso}
Nicolai Meinshausen and Bin Yu.
\newblock Lasso-type recovery of sparse representations for high-dimensional
  data.
\newblock {\em Annals of Statistics}, 37(1):246--270, 2009.

\bibitem{nardi}
Y.~Nardi and A.~Rinaldo.
\newblock On the asymptotic properties of the group lasso estimator for linear
  models.
\newblock {\em Electronic Journal of Statistics}, 2:605--633, 2008.

\bibitem{negahban}
Sahand~N. Negahban, Pradeep Ravikumar, Martin~J. Wainwright, and Bin Yu.
\newblock A unified framework for high-dimensional analysis of {$M$}-estimators
  with decomposable regularizers.
\newblock {\em Statist. Sci.}, 27(4):538--557, 2012.

\bibitem{nesterov1987}
Yurii Nesterov and Arkadii Nemirovskii.
\newblock {\em Interior-point polynomial algorithms in convex programming},
  volume~13 of {\em SIAM Studies in Applied Mathematics}.
\newblock Society for Industrial and Applied Mathematics (SIAM), Philadelphia,
  PA, 1994.

\bibitem{osborne2000}
M.~R. Osborne, Brett Presnell, and B.~A. Turlach.
\newblock A new approach to variable selection in least squares problems.
\newblock {\em IMA J. Numer. Anal.}, 20(3):389--403, 2000.

\bibitem{park}
Mee~Young Park and Trevor Hastie.
\newblock {$L_1$}-regularization path algorithm for generalized linear models.
\newblock {\em Journal of the Royal Statistical Society Series B},
  69(4):659--677, 2007.

\bibitem{ravikumar}
P.~Ravikumar, J.~Lafferty, H.~Liu, and L.~Wasserman.
\newblock Sparse additive models.
\newblock {\em Journal of the Royal Statistical Society Series B},
  71(5):1009--1030, 2009.

\bibitem{schwarz}
Gideon Schwarz.
\newblock Estimating the dimension of a model.
\newblock {\em Annals of Statistics}, 6(2):461--464, 1978.

\bibitem{TS}
Bernadetta Tarigan and Sara~A. van~de Geer.
\newblock Classifiers of support vector machine type with {$l_1$} complexity
  regularization.
\newblock {\em Bernoulli}, 12(6):1045--1076, 2006.

\bibitem{TR}
Robert Tibshirani.
\newblock Regression shrinkage and selection via the lasso.
\newblock {\em Journal of the Royal Statistical Society Series B},
  58(1):267--288, 1996.

\bibitem{van2008}
Sara~A. van~de Geer.
\newblock High-dimensional generalized linear models and the lasso.
\newblock {\em Annals of Statistics}, 36(2):614--645, 2008.

\bibitem{van2009}
Sara~A. van~de Geer and Peter B{\"u}hlmann.
\newblock On the conditions used to prove oracle results for the {L}asso.
\newblock {\em Electronic Journal of Statistics}, 3:1360--1392, 2009.

\bibitem{wu2009genome}
T.T. Wu, Y.F. Chen, T.~Hastie, E.~Sobel, and K.~Lange.
\newblock Genome-wide association analysis by lasso penalized logistic
  regression.
\newblock {\em Bioinformatics}, 25(6):714--721, 2009.

\bibitem{yuan2006}
Ming Yuan and Yi~Lin.
\newblock Model selection and estimation in regression with grouped variables.
\newblock {\em Journal of the Royal Statistical Society Series B},
  68(1):49--67, 2006.

\bibitem{zhang2008sparsity}
Cun-Hui Zhang and Jian Huang.
\newblock The sparsity and bias of the {LASSO} selection in high-dimensional
  linear regression.
\newblock {\em Annals of Statistics}, 36(4):1567--1594, 2008.

\bibitem{zhao2006}
Peng Zhao and Bin Yu.
\newblock On model selection consistency of {L}asso.
\newblock {\em J. Mach. Learn. Res.}, 7:2541--2563, 2006.

\bibitem{zou2006}
Hui Zou.
\newblock The adaptive lasso and its oracle properties.
\newblock {\em J. Amer. Statist. Assoc.}, 101(476):1418--1429, 2006.

\end{thebibliography}
\def\cprime{$'$}

\end{document}